\documentclass[12pt,a4paper,reqno]{amsart}
\usepackage{amssymb}
\usepackage{latexsym}
\usepackage{exscale}

\numberwithin{equation}{section}
\newtheorem{theor}{Theorem}[section]
\newtheorem{lemma}[theor]{Lemma}

\newtheorem{corol}[theor]{Corollary}
\newtheorem{remark}[theor]{Remark}
\newtheorem{prop}[theor]{Proposition}

\newcommand{\re}{\mathbb{R}}
\newcommand{\co}{\mathbb{C}}

\newcommand{\RR}{\mathbb{R}}
\newcommand{\NN}{\mathbb{N}}

\newcommand{\F}{\mathcal{F}}
\newcommand{\R}{\mathcal{R}}
\newcommand{\W}{\mathcal{W}}

\newcommand{\D}{\mathcal{D}}
\newcommand{\E}{\mathcal{E}}

\newcommand{\cals}{\mathcal{S}}
\newcommand{\calX}{\mathcal{X}}

\newcommand{\X}{\mathcal{X}}
\newcommand{\A}{\mathcal{A}}

\newcommand{\ep}{\varepsilon}

\newcommand{\bigchi}{\mathop{\mathchoice%
{\mbox{\Large$\chi$}}{\mbox{\large$\chi$}}{\mbox{\normalsize$\chi$}}%
{\mbox{\small$\chi$}}}\nolimits}

\newcommand{\BMO}{{\rm BMO}}

\renewcommand{\emptyset}{\mbox{\rm \O}}

\def\essinf{\mathop{\rm ess\ inf}}

\def\supp{\mathop{\rm supp}}

\newcommand{\aver}[1]{-\hskip-0.46cm\int_{#1}}

\begin{document}
\allowdisplaybreaks

\title[Weighted norm inequalities and elliptic
operators]{Weighted norm inequalities,  off-diagonal estimates and
elliptic operators
\\[.2cm]
{\footnotesize Part I: General operator theory and weights}}

\author{Pascal Auscher}

\address{Pascal Auscher
\\
Universit\'e de Paris-Sud et CNRS UMR 8628
\\
91405 Orsay Cedex, France} \email{pascal.auscher@math.p-sud.fr}

\author{Jos\'e Mar{\'\i}a Martell}

\address{Jos\'e Mar{\'\i}a Martell
\\
Instituto de Matem\'aticas y F{\'\i}sica Fundamental
\\
Consejo Superior de Investigaciones Cient{\'\i}ficas
\\
C/ Serrano 123
\\
28006 Madrid, Spain}

\address{\null\vskip-.7cm and\vskip-.7cm\null}

\address{Departamento de Matem\'aticas \\ Universidad Aut\'onoma de Madrid \\
28049 Madrid, Spain } \email{chema.martell@uam.es}

\thanks{This work was partially supported by the European Union
(IHP Network ``Harmonic Analysis and Related Problems'' 2002-2006,
Contract HPRN-CT-2001-00273-HARP). The second author was also
supported by MEC ``Programa Ram\'on y Cajal, 2005'' and by MEC Grant
MTM2004-00678.}

\date{\today}
\subjclass[2000]{42B20, 42B25}

\keywords{Good-$\lambda$ inequalities, Calder\'on-Zygmund
decomposition, Muckenhoupt weights, vector-valued inequalities,
extrapolation, singular non-integral operators, commutators with
bounded mean oscillation functions}

\begin{abstract}
This is the first part of a series of four articles. In this
work, we are interested in weighted norm estimates. We put the
emphasis on two results of different nature: one is based on a
good-$\lambda$ inequality with two-parameters and the other
uses Calder\'on-Zygmund decomposition. These results apply well to
singular ``non-integral'' operators and their commutators with
bounded mean oscillation functions. Singular means that they are
of order 0, ``non-integral'' that they do not have an integral
representation by a kernel with size estimates, even rough, so
that they may not be bounded on all $L^p$ spaces for $1<p<\infty$.
Pointwise estimates are then replaced by appropriate localized
$L^p-L^q$ estimates.  We obtain weighted $L^p$ estimates for a
range of $p$ that is different from $(1,\infty)$ and isolate the
right class of weights. In particular, we prove an extrapolation
theorem ``\`a la Rubio de Francia"  for such a class and thus
vector-valued estimates.
\end{abstract}

\maketitle

\begin{quote}
{\footnotesize\tableofcontents}
\end{quote}

\section*{General introduction}

This is a general introduction for this article and the series \cite{AM2, AM3, AM4}.

Calder\'on-Zygmund operators have been thoroughly studied since the
50's. They are singular integral operators associated with a kernel
satisfying certain size and smoothness conditions. One first shows
that the operator in question is bounded on $L^{p_0}$ for some $p_0$:
typically, for $p_0=2$ with  spectral theory, Fourier transform or
even the powerful $T(1)$, $T(b)$ theorems. Once this is achieved,
using the properties of the kernel,  one gets a weak-type (1,1)
estimate hence strong type $(p,p)$ for $1<p<p_{0}$  by means of the
Calder\'on-Zygmund decomposition and for $p>p_{0}$, one uses duality or
boundedness from $L^\infty$ to $\BMO$ and interpolation.  Still
another way for $p>p_{0}$ relies on
 good-$\lambda$
estimates via the Fefferman-Stein sharp maximal function. It is interesting to note that both Calder\'on-Zygmund decomposition and good-$\lambda$ arguments use independent  smoothness conditions on
the kernel, allowing generalizations in various ways.

The removal of regularity assumptions on the kernel is important, for
instance towards  applications to operators on non-smooth domains.
Let us mention  \cite{DMc} where a weak-type (1,1) criterion is
obtained under upper bound assumption on the kernel but no regularity in
the classical H\"older sense or in the sense of the H\"ormander condition
\cite{Hor}.

We mention also that Calder\'on-Zygmund operators satisfy also
commutator estimates with bounded mean oscillation functions and it
is therefore natural to try to extend them (see also the work of
\cite{DY1} in this direction following the methods in
\cite{DMc}).

A natural question is in what sense one should use the kernel of the
operators. It has become common practice but is it a necessary
limitation or a technical one. Indeed, one encounters
Calder\'on-Zygmund like operators  without any (reasonable) information
on their kernels which we call, following the implicit terminology
introduced in \cite{BK1}, singular ``non-integral''  operators in the
sense that they are still of order 0  but they do not have an
integral representation by a kernel with size and/or smoothness
estimates. The goal is to obtain some range of exponents $p$ for
which $L^p$ boundedness holds, and because this range may not be
$(1,\infty)$, one should abandon any use of kernels.

The first step was done  in  \cite{BK1}   where a criterion for weak-type $(p,p)$  for  some $p<p_{0}$ is presented.  In fact, this criterion is in the air in \cite{Fef}  but, still, \cite{BK1} brings some novelty such as the removal of the mean value property already observed  in \cite{DMc} when $p=1$. See also
\cite{BK2} and \cite{HM} for $L^p$ bounds $p<2$ of  the Riesz transforms of elliptic operators
starting from the $L^2$ bound proved in
\cite{AHLMcT}.

The second step was taken in \cite{ACDH}, inspired by the good-$\lambda$
estimates in the Ph.D.
thesis of one of us \cite{Ma1, Ma2}, where a criterion for strong
type $(p,p)$ for \textit{some} $p>p_{0}$   is proved and applied to Riesz transforms for the
Laplace-Beltrami operators on some Riemannian manifolds. A
criterion in the same spirit for a limited range of $p$'s also appears implicitly in
\cite{CP} towards perturbation theory for linear and non-linear elliptic
equations and  more explicitly  in \cite{Shen1, Shen2} (actually, we
shall observe here that the criterion in \cite{Shen2} is a corollary
of the one in \cite{ACDH}).

These two criteria are exposed in \cite{Aus-mem}, to  which the
reader is referred,  in the Euclidean setting  and applied to other
operators.

Our purpose is to investigate the weighted norm counterparts of this
new theory for Muckenhoupt weights and to apply this in the
subsequent papers. Again, the weighted norm theory is well known for
Calder\'on-Zygmund operators and we seek for criteria  applying to
larger classes of operators \textit{without kernel bounds} hence with
limited range of exponents. We mention \cite{Ma1} where some weighted
estimates for a functional calculi are proved but again assuming
appropriate kernel upper bounds. Our study will also clarify some
points in the unweighted case: in particular, we  present a simple
machinery  to prove (new) commutator estimates (both unweighted and
weighted) in this generality.

This paper is  concerned with the general operator theory and weights in the setting of spaces of homogeneous type. We study weighted boundedness criteria for operators and theirs commutators with bounded mean oscillation functions. Available machinery give us also vector-valued estimates. See the specific introductions of Parts \ref{part:one} and \ref{part:two} in this paper.

Part II, \cite{AM2},  is of independent interest as it develops a
theory of off-diagonal estimates in the context of spaces of
homogeneous type. In particular, the case of the semigroups generated
by elliptic operators is thoroughly studied. This is instrumental in
the application of the general theory in
\cite{AM3}.

In  Part III, \cite{AM3}, we consider operators arising from second
order elliptic operators $L$: operators of the type $\varphi(L)$ from
holomorphic functional calculus,  the Riesz transforms,  square
functions, \dots. We  obtain sharp or nearly sharp ranges of weighted
boundedness  of such operators, of  their commutators with  bounded mean oscillation functions,  and also
vector-valued inequalities.

In Part IV, \cite{AM4}, we apply our general theory to the Riesz
transform  on some Riemannian
manifolds or Lie groups as in \cite{ACDH} and their commutators.

\part{Good-$\lambda$ methods}\label{part:one}

\section{Introduction}

Good-$\lambda$ inequalities,  brought  to  Harmonic Analysis in
\cite{BG},   provide a powerful tool to prove
boundedness results for operators or at least comparisons of two
operators.  A typical good-$\lambda$ inequality  for two non-negative
functions $F$ and $G$ is as follows:  for every $0<\delta<1$ there
exists $\gamma=\gamma(\delta)$ and for every $w\in A_\infty$, there
exists  $0<\epsilon_w\le 1$ and $C_w>0$ such that for any $\lambda>0$
\begin{equation}\label{good-lambda:w:intro}
w\{x:F(x)>2\,\lambda, G(x)\le \gamma\,\lambda\}
\le
C_w\,\delta^{\epsilon_w}\, w\{x:F(x)>\lambda\}.
\end{equation}
The usual approach for proving such an estimate consists in  first
deriving a local version of it with respect to the underlying doubling
measure, and then passing to the weighted measure using that $w\in
A_\infty$.

Weighted good-$\lambda$ estimates encode a lot of information
about $F$ and $G$, since they give a comparison of the $w$-measure
of the level sets of both functions. As a consequence of
\eqref{good-lambda:w:intro} one gets, for instance, that for every
$0<p<\infty$ and all $w\in A_\infty$ then $\|F\|_{L^p(w)}$ is
controlled by $\|G\|_{L^p(w)}$. The same inequality holds with
$L^{p,\infty}$ in place of $L^p$ or with some other  function
spaces. Thus, the size of $F$ is controlled by that of $G$.

 In
applications, one tries to control a specific operator $T$ to be studied by a maximal one $M$  whose properties are known by setting
$F=Tf$ and $G=Mf$. For example,
  a Calder\'on-Zygmund operator  by
the Hardy-Littlewood maximal operator \cite{Coi}, \cite{CF}; a
fractional integral by  a fractional maximal operator \cite{MW}; a
Littlewood-Paley square function  by   a non-tangential maximal
operator  \cite{CWW}, \cite{Dah}, \cite{DJK},
\cite{GW}, \cite{Wil};  the maximal operator by  the sharp maximal
operator  \cite{FS}.

When $T$ is a Calder\'on-Zygmund operator with smooth kernel, in
particular it is already bounded on (unweighted) $L^2$, it was shown
in \cite{Coi}, \cite{CF} that \eqref{good-lambda:w:intro} holds with
$F=Tf$ and $G=M f$ with $M$ being the Hardy-Littlewood maximal
function. Thus,  $T$ is \lq\lq controlled\rq\rq\ by $M$ in $L^p(w)$
for all $0<p<\infty$ and $w\in A_{\infty}$ and therefore $T$ is
bounded on $L^p(w)$  if $M$ is bounded on $L^p(w)$, which by
Muckenhoupt's theorem means $w\in A_{p}$. In particular,  the range
of unweighted $L^p$ boundedness of $T$, that is the set  of $p$ for
which $T$ is  strong-type $(p,p)$,  is $(1,\infty)$, a fact that was
known by Calder\'on-Zygmund methods (see Part 2 of this paper).

Replacing $Mf$  by
$M(|f|^{p_0})^{1/p_0}$ for some $p_0>1$   changes the range of
unweighted $L^p$ boundedness   to $(p_{0}, \infty)$. See for instance
\cite{MPT}, and the references therein, where this occurs for
Calder\'on-Zygmund operators with less regular kernels.
  In this case, weighted $L^p(w)$ boundedness holds if
$w\in A_{p/p_{0}}$.

So far, there is a lower limitation on $p$ but no
upper limitation  in the sense that $p$ goes all the way to $\infty$.
This has to be so by a special and very simple case of Rubio de
Francia's extrapolation theorem (see \cite{Rub}, \cite{Gar}) which
says that any sublinear operator $T$ that is bounded on
$L^{p_{1}}(w)$ for some $0<p_1<\infty$ and all $w\in A_1$, is bounded
on $L^p$ for all $p_1\le p<\infty$.

Obviously, the above good-$\lambda$ inequality does not apply to
operators whose $L^p$ boundedness is expected  for $p_{0}<p<q_{0}$
with a finite exponent $q_{0}$. An example is  the Riesz transform
for the Laplace-Beltrami operator on some Riemannian manifolds
studied in \cite{ACDH, AC}. There,  a two-parameter good-$\lambda$
estimate incorporating an upper limitation in $p$ is used for proving
$L^p$ boundedness with a limited range of $p>2$.     See also
\cite{Aus-mem},
\cite{CP},  \cite{Shen1, Shen2}. These
two-parameter good-$\lambda$ estimates are of the form
\begin{equation}\label{good-lambda:two-param:intro}
\big|\{x: M F(x)> K\,\lambda, G(x)\le \gamma\,\lambda\}\big|
 \le
C\,
\left(
\frac{1}{K^{q_0}}+\frac{\gamma}{K}\right)\,
\big|\{x:M F(x)> \lambda\}\big|,
\end{equation}
for all $\lambda>0$, $K\ge K_0$ and $0<\gamma<1$. Note the explicit
dependance on $K, \gamma$ which are the two parameters and the
appearance of the exponent $q_{0} \in (0, \infty]$ in the right hand
side. From this, it follows  that $MF$ is controlled by $G$ in $L^p$
for all $0<p<q_0$.

The aim of this part  is to state conditions to obtain a weighted analog of  \eqref{good-lambda:two-param:intro}  and to derive some consequences for the study of operators. As we see below (Section \ref{sec:main}),
this  forces us to
specify  the power $\epsilon_{w}$  in
\eqref{good-lambda:w:intro}, hence to specify the reverse H\"older class for $w$. Indeed,  taking $w\in RH_{s'}$
then $\epsilon_{w}=1/s$  and we obtain the control of $MF$ by $G$ in
$L^p(w)$ for all $0<p<q_0/s$ (note that this implies that $w\in
RH_{(q_{0}/p)'}$).

This allows us to formulate simple unweighted conditions for the
$L^p(w)$ boundedness of (singular ``non-integral'') operators a
priori bounded on (unweighted) $L^p$ for $p_{0}<p<q_{0}$ for weights in the class
$\W^p(p_{0}, q_{0})= A_{p/p_{0}}\cap RH_{(q_{0}/p)'}$ (Section \ref{section:SIO}). A slight
improvement furnishes, almost for free, boundedness of their
commutators with bounded mean oscillation functions for the same
weights (Section \ref{section:comm:p-big}). This class of weights (studied in Section \ref{subsection:set:W}) is the largest possible within
$A_{\infty}$ as we prove an extrapolation result for it. Namely, if $T$ is
bounded on some $L^p(w)$ for some fixed $p$ and for all $w\in
\W^p\big(p_{0}, q_{0}\big)$, then the same happens for every $q\in
(p_0,q_0)$ and the corresponding class of weights. Using
  ideas on extrapolation from \cite{CMP} and \cite{CGMP},  we obtain
vector-valued inequalities automatically again for limited ranges of
$p$ (Section \ref{subsection:extrapo}).  For simplicity of the exposition, we work in
the Euclidean space equipped with the Lebesgue measure. See Section
\ref{section:SHT}  for extensions to spaces of homogeneous type.

\section{Muckenhoupt weights}\label{subsection:weights}

We review some needed background  on Muckenhoupt weights. We
use the notation
$$
\aver{E} h = \frac {1}{|E|} \int_{E} h(x)\, dx
$$
and we often forget  the Lebesgue measure  and the variable of the integrand in writing integrals, unless this is needed to avoid confusions.

A weight $w$ is a non-negative locally integrable function. We say
that $w\in A_p$, $1<p<\infty$, if there exists a constant $C$ such
that for every ball $B\subset\re^n$ (balls could be switched to cubes)
$$
\Big(\aver{B} w\Big)\,
\Big(\aver{B} w^{1-p'}\Big)^{p-1}\le C.
$$
For $p=1$, we say that $w\in A_1$ if there is a constant $C$ such
that for every ball $B\subset \re^n$
$$
\aver{B} w
\le
C\, w(x),
\qquad \mbox{for a.e. }x\in B,
$$
or, equivalently, $M w\le C\,w$  a.e. where $M$ denotes the
uncentered maximal operator over balls (or cubes) in $\re^n$. The
reverse H\"older classes are defined in the following way: $w\in
RH_{q}$, $1< q<\infty$, if there is a constant $C$ such that for
every ball $B\subset \re^n$
$$
\Big(\aver{B} w^q\Big)^{\frac1q}
\le C\, \aver{B} w.
$$
The endpoint $q=\infty$ is given by the condition: $w\in RH_{\infty}$
whenever, for any ball $B$,
$$
w(x)\le C\, \aver{B} w,
\qquad \mbox{for a.e. }x\in B.
$$
Notice that we have excluded the case $q=1$ since the class $RH_1$
consists of all the weights, and that is the way $RH_1$ is
understood in what follows.

We sum up some of the properties of these classes in the following result.
\begin{prop}\label{prop:weights}\
\begin{enumerate}
\renewcommand{\theenumi}{\roman{enumi}}
\renewcommand{\labelenumi}{$(\theenumi)$}
\addtolength{\itemsep}{0.2cm}

\item $A_1\subset A_p\subset A_q$ for $1\le p\le q<\infty$.

\item $RH_{\infty}\subset RH_q\subset RH_p$ for $1<p\le q\le \infty$.

\item If $w\in A_p$, $1<p<\infty$, then there exists $1<q<p$ such
that $w\in A_q$.

\item If $w\in RH_q$, $1<q<\infty$, then there exists $q<p<\infty$ such
that $w\in RH_p$.

\item $\displaystyle A_\infty=\bigcup_{1\le p<\infty} A_p=\bigcup_{1<q\le
\infty} RH_q . $

\item If $1<p<\infty$, $w\in A_p$ if and only if $w^{1-p'}\in
A_{p'}$.

\item If $1\le q\le \infty$ and $1\le s<\infty$, then $\displaystyle
w\in A_q \cap RH_s$ if and only if $ w^{s}\in A_{s\,(q-1)+1}$.

\end{enumerate}
\end{prop}

Properties $(i)$-$(vi)$ are standard, see for instance \cite{GR} or
\cite{Duo}. For $(vii)$ see \cite{JN}.

\section{Two parameter good-$\lambda$ estimates}\label{section:good-lambda}

Unless specified otherwise, $M$ denotes the uncentered maximal operator over cubes (or balls) in $\re^n$.

\subsection{Main result}\label{sec:main}

\begin{theor}\label{theor:good-lambda:w}
Fix $1<q\le \infty$, $a\ge 1$ and $w\in RH_{s'}$, $1\le s<\infty$.
Then, there exist $C=C(q,n,a,w,s)$ and $K_0=K_0(n,a)\ge 1$ with the
following property: Assume that $F$, $G$, $H_1$ and $H_2$ are
non-negative measurable functions  on $\re^n$ such that for any cube $Q$ there
exist non-negative functions $G_Q$ and $H_Q$ with $F(x)\le G_Q(x)+
H_Q(x)$ for a.e. $x\in Q$ and
\begin{equation}\label{H-Q}
\Big(\aver{Q} H_Q^q\Big)^{\frac1q}
\le
a\, \big(M F(x) + M H_1(x) + H_2(\bar{x})\big),
\qquad
\forall\,x,\bar{x}\in Q;
\end{equation}
and
\begin{equation}\label{G-Q}
\aver{Q} G_Q
\le
G(x),
\qquad
\forall\,x\in Q.
\end{equation}
Then for all $\lambda>0$, $K\ge K_0$ and $0<\gamma<1$
\begin{equation}
w
\big\{M F> K\,\lambda, G+H_2\le \gamma\,\lambda\big\}
 \le
C\,
\left(
\frac{a^q}{K^q}+\frac{\gamma}{K}\right)^{\frac1s} \, w
\big\{M F+ M H_1> \lambda\big\}.
\label{good-lambda:w}
\end{equation}
As a consequence, for all $0<p
<\frac{q}{s}$, we have
\begin{equation}\label{good-lambda:Lp:w}
\|M F\|_{L^p(w)}
\le
C\,\big(\|G\|_{L^p(w)} + \|M H_1\|_{L^p(w)}+\|H_2\|_{L^p(w)}\big),
\end{equation}
provided $\|M F\|_{L^p(w)}<\infty$, and
\begin{equation}\label{good-lambda:Lp-weak:w}
\|M F\|_{L^{p,\infty}(w)}
\le
C\,\big(\|G\|_{L^{p,\infty}(w)} + \|M
H_1\|_{L^{p,\infty}(w)}+\|H_2\|_{L^{p,\infty}(w)}\big),
\end{equation}
provided $\|M F\|_{L^{p,\infty}(w)}<\infty$. Furthermore, if $p\ge
1$ then \eqref{good-lambda:Lp:w} and \eqref{good-lambda:Lp-weak:w}
hold, provided $F\in L^1$ \textup{(}whether or not $M F\in
L^p(w)$\textup{)}.
\end{theor}

The proof of this result is in Section
\ref{subs:proof:theor:good-lambda:w}.

\begin{remark}\rm
We do mean that the estimates \eqref{H-Q} and \eqref{G-Q} are
valid at \emph{any} points $x, \bar x \in Q$, not just almost
everywhere.
\end{remark}

\begin{remark}\label{remark:q=infty}\rm
The case $q=\infty$ is the standard one:  the $L^q$-average appearing
in the hypothesis is understood as an essential supremum and
$K^{-q}=0$. Thus,  the $L^p(w)$ and $L^{p,\infty}(w)$ estimates will
hold for any $0<p<\infty$, no matter the value of $s$, that is, for
any $w\in A_\infty$.
\end{remark}

\begin{remark}\label{remark:q-large}\rm
If \eqref{H-Q} holds for any $q>1$, then \eqref{good-lambda:Lp:w}
holds for all $0<p<\infty$ and for all $w\in A_\infty$. To see
this, we fix $0<p<\infty$ and $w\in A_\infty$. Then $w\in RH_{s'}$
for some $1\le s<\infty$ and it suffices to take $q$ large enough
so that $p<q/s$.
\end{remark}

\begin{remark}\rm
In applications, error terms appear in localization arguments either
in the form $M H_1(x)$ or $H_2(\bar{x})$ (with $\bar{x}$ independent
of $x$) or both. The unweighted case
\cite[Theorem 2.4]{ACDH} is of this type.
\end{remark}

\begin{remark}\label{remark:up-endpoint}\rm
If $s>1$ and $q<\infty$, then one also obtains the end-point $p=q/s$.
To do it, we only need to observe that $w\in RH_{s_0'}$ for some
$1<s_0<s$ (see $(v)$ in Proposition \ref{prop:weights}) and so we can
apply Theorem \ref{theor:good-lambda:w} with $p=q/s<q/s_0$.
\end{remark}

We present some  applications of Theorem
\ref{theor:good-lambda:w} recovering some previously  known
estimates.

\subsection{Fefferman-Stein Inequality}
The classical Fefferman-Stein inequality relating  $M$ and $M^\#$
follows at once from Theorem \ref{theor:good-lambda:w}. We take
$F=|f|\in L^{1}_{\rm loc}$, $H_1=H_2=0$. For each cube $Q$ we denote
by $f_{Q}$ the average of $f$ on $Q$,
$$
F= |f|
\le
|f_Q|+|f-f_Q|\equiv H_Q+G_Q.
$$
Taking $q=\infty$, we trivially have
$\|H_Q\|_{L^{\infty}(Q)}=|f_Q|\le Mf(x)= MF(x)$ for each $x\in Q$. Also, by definition of $M^\#$
$$
\aver{Q} G_Q
=
\aver{Q} |f-f_Q|
\le M^\# f(x)\equiv G(x),
\qquad
\forall\,x\in Q.
$$
Thus, \eqref{good-lambda:w} holds (with $q=\infty$) and consequently,
for every $0<p<\infty$ and every $w\in A_\infty$ we have
\begin{equation}\label{Fef-Ste}
\|M f\|_{L^p(w)}
\le
C\,\|M^\# f\|_{L^p(w)},
\end{equation}
whenever $Mf\in L^p(w)$. This is what is proved in \cite{FS}.

\subsection{Generalized sharp maximal functions}

In \cite{Ma1}, a generalization of $M^\#$ is introduced in the
setting of spaces of homogeneous type. In the Euclidean setting, we
define $M^{\#}_D$ as follows. Let $\{D_t\}_{t>0}$ be a family of
operators (for instance, an approximation of the identity but it
could be more general) such that each $D_t$ is an integral operator
with kernel $d_t(x,y)$ for which
$$
|d_t(x,y)|\le C\,t^{-\frac{n}{m}}\,h\big(|x-y|^m\,t^{-1}\big)
$$
where $m$ is some positive fixed constant and $h$ is positive,
bounded, decreasing and decaying to $0$ fast enough. Then we define a
new sharp maximal function associated to $\{D_t\}_{t>0}$ as
$$
M^{\#}_D f(x)
=
\sup_{Q\ni x}\aver{Q} |f-D_{t_Q} f|
$$
where $t_Q=\ell(Q)^m$ and $\ell(Q)$ is the sidelength of $Q$.

Examples  are given by the semigroups associated with a second order
elliptic operators $\{e^{-t\,L}\}_{t>0}$ whose heat kernels have
Gaussian (or some other) decay (see \cite{AT, Ma1, DMc, ArTe}, \dots)

With Theorem \ref{theor:good-lambda:w} we can reprove the
good-$\lambda$ inequality of \cite{Ma1} for $M^\#_D$ and $M$. As
before take $F=|f|\in L^{p}$ for some $p\ge 1$, $H_1=H_2=0$. For each
cube $Q$ we write
$$
F= |f|
\le
|D_{t_Q}f|+|f-D_{t_Q}f|\equiv H_Q+G_Q.
$$
Taking $q=\infty$, we  have $\|H_Q\|_{L^{\infty}(Q)}\le C\, Mf(x)=
C\, MF(x)$ for each $x\in Q$ by the properties assumed on $D_{t}$. Moreover, by definition of
$M^\#_D$,
$$
\aver{Q} G_Q
=
\aver{Q} |f-D_{t_Q}f|
\le M^\#_D f(x)\equiv G(x),
\qquad
\forall\,x\in Q.
$$
Thus, one obtains \eqref{good-lambda:w} (with $q=\infty$) and hence,
for every $0<p<\infty$ and every $w\in A_\infty$ we have
$$
\|M f\|_{L^p(w)}
\le
C\,\|M^\#_D f\|_{L^p(w)},
$$
whenever $Mf\in L^p(w)$. This is the result proved in \cite{Ma1}.

\subsection{Applications to Singular ``Non-Integral'' Operators}\label{section:SIO}

We present here different applications of Theorem
\ref{theor:good-lambda:w} toward   weighted norm inequalities for
operators,  avoiding all use of kernel representation, hence the
terminology ``non-integral".

In what follows, we say that an operator  $T$ acts from $A$ into $B$
(with $A$, $B$ being some given sets) if $T$ is a map defined on $A$
and valued in $B$. An operator $T$ acting from $A$ to $B$, both
vector spaces of measurable functions,  is sublinear if $$\vert
T(f+g)\vert \le \vert Tf\vert + \vert Tg\vert  \qquad {\rm and}\qquad
\vert T(\lambda f)\vert = \vert \lambda\vert \, \vert Tf\vert$$ for
all $f,g \in A$ and $\lambda\in \re$ or $\co$. Let us mention that
for the theorems of this section, the second condition is  not
needed.

\begin{theor} \label{theor:main-w}
Let $1\le p_0<q_0\le \infty$. Let $\E$ and $\D$ be vector spaces
such that $\D\subset \E$.  Let $T$, $S$ be  operators such that
$S$ acts from $\D$ into the set of measurable functions and $T$ is
sublinear acting from $\E$ into $L^{p_0}$. Let $\{\A_r\}_{r>0}$ be
a family of operators acting from $\D$ into $\E$. Assume that
\begin{equation}\label{T:I-A}
\Big(\aver{B} |T(I-\A_{r(B)})f|^{p_0}\Big)^{\frac1{p_0}}
\le
C\, M\big(\,|S f|^{p_0}\big)^{\frac1{p_0}}(x),
\end{equation}
and
\begin{equation}\label{T:A}
\Big(\aver{B} |T\A_{r(B)}f|^{q_0}\Big)^{\frac1{q_0}}
\le
C\, M\big(\,|T f|^{p_0}\big)^{\frac1{p_0}}(x),
\end{equation}
for all $f\in \D$, all ball $B$ where $r(B)$ denotes its radius and
all $x\in B$. Let $p_0<p<q_0$ \textup{(}or $p=q_0$ when
$q_0<\infty$\textup{)} and $w\in A_{\frac{p}{p_0}}\cap
RH_{\left(\frac{q_0}{p}\right)'}$. There is a constant $C$ such that
\begin{equation}\label{est:main-th}
\|T f\|_{L^p(w)}
\le
C\, \|S f\|_{L^p(w)}
\end{equation}
for all $f\in \D$. Furthermore,
  for all $p_0<r<q_0$, there is a constant $C$ such that
\begin{equation}\label{T:v-v}
\Big\|
\Big(\sum_j |T f_j|^r\Big)^{\frac1r}
\Big\|_{L^p(w)}
\le
C\,\Big\|
\Big(\sum_j |S f_j|^r\Big)^{\frac1r}
\Big\|_{L^p(w)}
\end{equation}
for all $f_j\in \D$.
\end{theor}

We would like to emphasize that  \eqref{T:I-A} and \eqref{T:A} are
unweighted assumptions. This is a triple extension of  \cite[Theorem
3.1]{ACDH}: we introduce a second operator $S$, obtain weighted
inequalities and also vector-valued estimates.

\begin{remark}\label{remark:S}\rm
The most common situation is $S=I$, $\E=L^{p_{0}}$ with $\D$ being a
class of \lq\lq nice\rq\rq\ functions such as $L^{p_0}$, $L^{p_0}\cap
L^2$, $L^{\infty}_c$, $C_0^\infty$, \dots. In that case,
\eqref{est:main-th}  is interesting only  when the right hand side is
finite, hence we may also impose  $f\in L^p(w)$. This implies the
boundedness of $T$ from $\D\cap L^p(w)$ into $L^p(w)$ for the
$L^p(w)$ norm. See \cite{AM3} for a situation where $S\ne I$.
\end{remark}

\begin{remark}\label{remark:q0=inf}\rm
In this result, the case $q_0=\infty$ is understood in the sense that
the $L^{q_0}$-average in \eqref{T:A} is indeed an essential supremum.
Besides, the condition for the weight turns out to be $w\in
A_{p/p_0}$ for $p>p_0$. Similarly, if
\eqref{T:A} is satisfied for all $q_0<\infty$ then
\eqref{est:main-th} holds for all $p_0<p<\infty$ and for all $w\in
A_{p/p_0}$.
\end{remark}

\begin{remark}\label{remark:A(B)}\rm
 A slightly more general statement consists in replacing the
family $\{\A_{r}\}$ by $\{\A_{B}\}$ indexed by balls. We use this below.\end{remark}

\begin{proof}[Proof of Theorem \ref{theor:main-w}]
The vector-valued inequalities \eqref{T:v-v} follow automatically by
extrapolation, see Theorem \ref{theor:extrapol} below.

We  prove \eqref{est:main-th},   first  in the case $q_0<\infty$ and
$p_0<p\le q_0$. Let $f\in \D$ and so $F=|T f|^{p_0}\in L^{1}$. Fix a
cube $Q$ (we switch to cubes for the proof). As $T$ is sublinear, we
have
$$
F
\le
G_Q+ H_Q \equiv 2^{p_0-1}\,|T(I-\A_{r(Q)})f|^{p_0} +
2^{p_0-1}\,|T\A_{r(Q)}f|^{p_0}.
$$
Then \eqref{T:I-A} and \eqref{T:A} yield the corresponding conditions
\eqref{H-Q} and \eqref{G-Q} with   $q=q_0/p_0$, $H_1=H_2\equiv 0$,
$a=2^{p_0-1}\,C^{p_0}$ and   $G=2^{p_0-1}\,C^{p_0}\, M
\big(\,|S f|^{p_0}\big)$. As $w\in RH_{(q_0/p)'}$,
Theorem \ref{theor:good-lambda:w} and Remark
\ref{remark:up-endpoint} (since $q_0<\infty$ implies $q<\infty$)
with $p/p_0>1$ in place of $p$ and $s=q_0/p$ yield
$$
\|T f\|_{L^{p}(w)}^{p_0}
\le
\|M F\|_{L^{\frac{p}{p_0}}(w)}
\le
C\,\|G\|_{L^{\frac{p}{p_0}}(w)}
=
C\,\big\|M
\big(\,|S f|^{p_0}\big)\big\|_{L^{\frac{p}{p_0}}(w)}
\le
C\,\|S f\|_{L^{p}(w)}^{p_0},
$$
where in the last estimate we have used that $w\in A_{p/p_0}$.

In the case $q_0=\infty$ and $p<\infty$, Theorem
\ref{theor:good-lambda:w} applies as before when $w\in A_{p/p_{0}}$
by  Remark \ref{remark:q=infty}.
\end{proof}

\begin{remark}\label{remark:weakt-type}\rm
Under the assumptions of Theorem \ref{theor:main-w}, we can also
prove an end-point weak-type estimate. Namely, if $w\in A_{1}\cap
RH_{\left(\frac{q_0}{p_0}\right)'}$, then there is a constant $C$
such that
\begin{equation}\label{est:main-th:weak}
\|T f\|_{L^{p_0,\infty}(w)}
\le
C\, \|S f\|_{L^{p_0}(w)},
\end{equation}
for all $f\in \D$. The proof follows the same ideas but one has to
use the weak-type estimate \eqref{good-lambda:Lp-weak:w} in place of
\eqref{good-lambda:Lp:w}. The details are left to the reader.

Let us recall that we have assumed that for $f\in \D$ then $F=|T
f|^{p_0}\in L^1$. This hypothesis is not granted directly for $T$ in
some applications (for instance, it is not true for $p_0=1$ and $T$
being the Hilbert transform or the Riesz transforms) but for suitable
approximations $T_\varepsilon$ that are bounded on $L^{p_0}(w)$ (with
some bound that is allowed to depend on $\varepsilon$).  In such a
case,  one obtains the weak-type estimate for $T_\varepsilon$ with a
uniform control on the constant and the weak-type estimate for $T$
follows by a limiting procedure. (This happens for the Hilbert
transform: the kernel is truncated in such a way that it is in $L^1$,
so the approximations $T_\varepsilon$ are bounded on $L^1$.) Let us
mention that for Calder\'on-Zygmund operators the usual approach is
different: the weighted weak-type $(1,1)$ estimate for $A_1$ weights
follows by using the Calder\'on-Zygmund decomposition (see
\cite[Chapter IV]{GR}), see also \cite{BK1} for a weak-type
$(p_0,p_0)$ with $p_0>1$, and Theorems \ref{theor:B-K:small},
\ref{theor:B-K:small:w} below.
\end{remark}

\begin{remark}\rm
Theorem \ref{theor:good-lambda:w} implies a variant of Theorem
\ref{theor:main-w} valid for all $0<p_0<q_0\le\infty$. We do not
know, however, whether such a result is useful in applications when
$p_0<1$. The precise statement and the minor modifications in the
proof are left to the reader.
\end{remark}

The following extension of Theorem \ref{theor:main-w} is also
useful. For simplicity we assume that $S=I$.

\begin{theor} \label{theor:ACDH:w}
Let $1\le p_0<q_0\le \infty$. Let $\D$, $\E$, $T$ and
$\{\A_r\}_{r>0}$ be as in Theorem \ref{theor:main-w}. Assume that
\eqref{T:I-A} holds with $S=I$ and, in place of \eqref{T:A},  that
\begin{equation}\label{T:A:ACDH}
\Big(\aver{B} |T\A_{r(B)}f|^{q_0}\Big)^{\frac1{q_0}}
\le
C\, \big(M\big(\,|T f|^{p_0}\big)^{\frac1{p_0}}(x) + M\big(\,|S_1
f|^{p_0}\big)^{\frac1{p_0}}(x) + |S_2 f(\bar{x})|\,\big),
\end{equation}
holds for all $f\in \D$ and  all $x$, $\bar{x}\in B$ where $S_1$,
$S_2$ are two given operators. Let $p_0<p<q_0$ and $w\in
A_{\frac{p}{p_0}}\cap RH_{\left(\frac{q_0}{p}\right)'}$. If $S_1$ and
$S_2$ are bounded on $L^{p}(w)$, then $$ \|T f\|_{L^p(w)}
\le C\, \|f\|_{L^p(w)}
$$
for all $f\in \D\cap L^p(w)$.
\end{theor}

Observe that Remarks \ref{remark:q0=inf} and \ref{remark:A(B)} apply
to this result. Also, the operator $T$ satisfies the vector-valued
inequalities \eqref{T:v-v}.

\begin{proof}
The proof is almost identical to the one of Theorem
\ref{theor:main-w}. Let $f\in \D\cap L^p(w)$ and set $F=
|Tf|^{p_{0}}\in L^1$, $H_1= |S_1 f|^{p_0}$ and $H_2=|S_2
f|^{p_0}$. Theorem \ref{theor:good-lambda:w} gives us
\begin{align*}
\|T f\|_{L^{p}(w)}^{p_0}
&
\le
\|M F\|_{L^{\frac{p}{p_0}}(w)}
\le
C\,\big(\|G\|_{L^{\frac{p}{p_0}}(w)} +\|M
H_1\|_{L^{\frac{p}{p_0}}(w)} + \|H_2\|_{L^{\frac{p}{p_0}}(w)}
\big)
\\
&=
C\,
\big(
\big\|M \big(\,|S f|^{p_0}\big)\big\|_{L^{\frac{p}{p_0}}(w)} +
\big\|M \big(\,|S_1 f|^{p_0}\big)\big\|_{L^{\frac{p}{p_0}}(w)}+
\big\||S_2 f|^{p_0}\big\|_{L^{\frac{p}{p_0}}(w)} \big)
\\
&
\le
C\,\|f\|_{L^{p}(w)}^{p_0},
\end{align*}
where we have used that $M$ is bounded on $L^{\frac{p}{p_0}}(w)$
(since $w\in A_{p/p_0}$) and that, by hypothesis, $S_1$, $S_2$ are
bounded on $L^p(w)$.
\end{proof}

The last result of this section is an extension of \cite[Theorem
3.1]{Shen2}.

\begin{theor} \label{theor:shen}
Let $1\le p_0<q_0\le \infty$. Suppose that $T$ is a bounded sublinear
operator  on $L^{p_0}$. Assume that there exist constants
$\alpha_{2}>\alpha_{1}>1$, $C>0$ such that
\begin{equation}\label{T:shen}
\Big(\aver{B} |Tf|^{q_0}\Big)^{\frac1{q_0}}
\le
C\, \bigg\{ \Big(\aver{\alpha_{1}\, B}
|Tf|^{p_0}\Big)^{\frac1{p_0}} +
M\big(|f|^{p_0}\big)^{\frac1{p_0}}(x)\bigg\},
\end{equation}
for all balls $B$, $x\in B$ and  all $f\in L^{\infty}$ with compact
support in $\RR^n \setminus \alpha_{2}\, B$. Let $p_0<p<q_0$
and $w\in A_{\frac{p}{p_0}}\cap RH_{\left(\frac{q_0}{p}\right)'}$.
Then, there is a constant $C$ such that
$$
\|T f\|_{L^p(w)}
\le
C\, \|f\|_{L^p(w)}
$$
for all    $f\in L^\infty$ with compact support.
\end{theor}

\begin{proof}
For any ball $B$, let $\A_{B}f= (1-\bigchi_{\alpha_{2}\, B})\, f$. We
fix $f\in L^\infty_c$, a ball $B$ and $x\in B$. Using the
$L^{p_{0}}$-boundedness of $T$, we have
\begin{equation}
\label{eq:shen}
\Big(\aver{\alpha_{1}\, B} |T(I-\A_{B})f|^{p_0}\Big)^{\frac1{p_0}}
\le
C\,
\Big(\aver{\alpha_{2}\, B}
|f|^{p_0}\Big)^{\frac1{p_0}}
\le
C\, M\big(|f|^{p_0}\big)^{\frac1{p_0}}(x).
\end{equation}
In particular \eqref{T:I-A} holds since $\alpha_1>1$. Next,  by
\eqref{T:shen} and since $|\A_B f|\le |f|$ we have
$$
\Big(\aver{B}  |T\A_{B}f|^{q_0}\Big)^{\frac1{q_0}} \le
C\, \bigg\{ \Big(\aver{\alpha_{1}\, B}
|T\A_{B}f|^{p_0}\Big)^{\frac1{p_0}} +
M\big(|f|^{p_0}\big)^{\frac1{p_0}}(x)\bigg\}.
$$
By \eqref{eq:shen} and the sublinearity of $T$, we obtain
$$
\Big(\aver{ B} |T\A_{B}f|^{q_0}\Big)^{\frac1{q_0}} \le
C\, M\big(|Tf|^{p_0}\big)^{\frac1{p_0}}({x})+ C\,
M\big(|f|^{p_0}\big)^{\frac1{p_0}}(x),
$$
which is \eqref{T:A:ACDH} with $S_{1}=I$ and $S_{2}=0$. We conclude
on applying Theorem \ref{theor:ACDH:w} with $\D=L^\infty_c$ and
$\E=L^{p_0}$.
\end{proof}

\subsection{Commutators with BMO functions: part I}\label{section:comm:p-big}

A slight strengthening  of the hypotheses in Theorem
\ref{theor:main-w} furnishes weighted $L^p$ estimates for commutators
with $\BMO$ functions.

Let $b\in \BMO$ (BMO is for bounded mean oscillation), that is,
$$
\|b\|_{\BMO}
=
\sup_B \aver{B} |b(x)-b_B|\, dx
<\infty,
$$
where the supremum is taken over all balls and $b_B$ stands for the
average of $b$ on $B$.
Let $T$ be a sublinear bounded operator on some $L^{p_{0}}$.
Boundedness is assumed to avoid technical issues with the definition of the commutators. It   could be relaxed, for instance, by imposing that $T$
acts from $\E=\cap_p L^{p}_c$ into $L^{p_0}$. Sublinearity is defined in Section
\ref{section:SIO}.

For any $k\in \NN$ we define the $k$-th order commutator
$$
T_b^k f(x)=T\big((b(x)-b)^k\,f\big)(x),
\qquad f\in L^\infty_c, \qquad x\in \RR^n.
$$
Note that $T_b^0=T$. Commutators are usually considered for linear
operators $T$ in which case they can be alternatively defined by
recurrence: the first order commutator is
$$
T_b^1f(x)=[b,T]f(x)= b(x)\,T f(x)- T(b\,f)(x)
$$
and for $k\ge 2$,  the $k$-th order commutator is given by
$T_b^k=[b,T_b^{k-1}]$.

We claim that since $T$ is bounded in $L^{p_0}$ then $T_b^k f$
is well defined  in $
L^{q}_{\rm loc}$ for any $0<q<p_0$ and
 for any $f\in L^\infty_c$: take a cube $Q$ containing
the support of $f$ and observe that by sublinearity for a.e. $x\in \re^n$
\begin{align*}
|T_b^k f(x)|
&\le
\sum_{m=0}^k
C_{m,k}\,|b(x)-b_Q|^{k-m}\,\big|T\big((b-b_Q)^m\,f\big)(x)\big|.
\end{align*}
 John-Nirenberg's inequality implies
$$
\int_Q |b(y)-b_Q|^{m\,p_0}\,|f(y)|^{p_0}\,dy
\le C
\|f\|_{L^\infty}\,\|b\|_{\BMO}^{m\,p_0}\,|Q|<+\infty.
$$
Hence, $T\big((b-b_Q)^m\,f\big) \in L^{p_{0}}$ and the claim follows.

We are going to see that  Theorem \ref{theor:good-lambda:w} can be
applied to $T_b^1$ where the function $H_2$ involves $T=T_{b}^0$. The
same will be done for $T_b^k$ and in this case $H_2$ involves the
preceding commutators $T, T_b^1, \dots, T_b^{k-1}$. Thus an induction
argument (details are  in Section \ref{subsec:proof:theor:comm})
will lead us to the following estimates:

\begin{theor} \label{theor:comm}
Let $1\le p_0<q_0\le \infty$ and $k\in \mathbb{N}$. Suppose that $T$
is a sublinear operator bounded on $L^{p_0}$, and that
$\{\A_r\}_{r>0}$ is a family of operators acting from $L^{\infty}_c$
into $L^{p_0}$. Assume that
\begin{equation}\label{T:I-A:comm}
\Big(\aver{B} |T(I-\A_{r(B)})f|^{p_0}\Big)^{\frac1{p_0}}
\le
C\,\sum_{j=1}^\infty \alpha_j\,\Big(\aver{2^{j+1}\,B}
|f|^{p_0}\Big)^{\frac1{p_0}},
\end{equation}
and
\begin{equation}\label{T:A:comm}
\Big(\aver{B} |T\A_{r(B)}f|^{q_0}\Big)^{\frac1{q_0}}
\le
\sum_{j=1}^\infty \alpha_j\,\Big(\aver{2^{j+1}\,B}
|Tf|^{p_0}\Big)^{\frac1{p_0}},
\end{equation}
for all $f\in L^{\infty}_c$ and all ball $B$ where $r(B)$ denotes its
radius. Let $p_0<p<q_0$ and $w\in A_{\frac{p}{p_0}}\cap
RH_{\left(\frac{q_0}{p}\right)'}$. If $\sum_{j} \alpha_j\,j^k<\infty$
then there is a constant $C$  such that for all $f\in L^\infty_{c}$ and all $b\in \BMO$,
\begin{equation}\label{comm-Lp:w}
\|T_b^k f\|_{L^p(w)}
\le C\, \|b\|_{\BMO}^k\,\|f\|_{L^p(w)},
\end{equation}
for all $f\in L^\infty_c$.
\end{theor}

\begin{remark}\label{remark:indu-0}\rm
Under the  assumptions above,  we have  $\sum_j \alpha_j<\infty$ and
so \eqref{T:I-A:comm} and \eqref{T:A:comm} imply respectively
\eqref{T:I-A} and \eqref{T:A}. Consequently, Theorem
\ref{theor:main-w} applies to $T=T_b^0$ and yields its
$L^p(w)$-boundedness.
\end{remark}

Observe that Remarks \ref{remark:q0=inf} and \ref{remark:A(B)} apply
to this result. Also, the operator $T_b^k$ satisfies the
vector-valued inequalities \eqref{T:v-v}. The assumptions
\eqref{T:I-A:comm} and \eqref{T:A:comm} can be relaxed in the spirit
of   Theorem \ref{theor:ACDH:w} by allowing error terms in the right hand sides:  details and proof are left to
the interested reader.

\begin{remark}\label{remark:multi-comm:p-big}\rm
As in \cite{PT} one can linearize the $k$-th order commutator and
consider the following multilinear commutators
$$
T_{\vec{b}} f(x)
=
T\Big( \Big(\prod_{j=1}^k (b_j(x)-b_j)\Big)\,f\Big) (x).
$$
where $\vec{b}=\{b_1,\dots,b_k\}$ is a family of $\BMO$ functions.
Notice that if $b_1=\cdots=b_k=b$ we have that $T_{\vec{b}}=T_b^k$.
The proof of Theorem \ref{theor:comm} can be adapted to $T_{\vec{b}}$
and thus get the corresponding weighted estimates for it (see Remark
\ref{ref:proof::multi-comm:p-big}). The precise statement is left to
the reader.
\end{remark}

\section{The sets $\W_w(p_0,q_0)$ and Extrapolation}\label{section:properties:w}

\subsection{The sets $\W_w(p_0,q_0)$}\label{subsection:set:W}

The conclusion of Theorem \ref{theor:main-w} with $S=I$ and $\D=L^{p_{0}}$ (and also of Theorems
\ref{theor:ACDH:w} and \ref{theor:comm}) can be rewritten as follows:
given $w\in A_\infty$, we introduce the set
$$
\W_w(p_0,q_0)
=
\big\{
p: p_0<p<q_0, w\in A_{\frac{p}{p_0}}\cap
RH_{\left(\frac{q_0}{p}\right)'}
\big\},
$$
and we have shown that $T$ is bounded on  $L^p(w)$   whenever $p\in
\W_w(p_0,q_0)$. Let us give some properties of this set.

\begin{lemma}\label{lemma:interval-w}
Let $w\in A_\infty$ and $1\le p_0<q_0\le \infty$. Then
$\W_w(p_0,q_0)=\big(\,p_0\,r_w, \frac{q_0}{(s_w)'}\big)$ where
$$
r_w=\inf\{r\ge 1:w\in A_r\},
\qquad\qquad
s_w=\sup\{s>1: w\in RH_s\}.
$$
\end{lemma}

If $q_0=\infty$, this result has to be understood in the following
way: the set $\W_w(p_0,q_0)$ is defined by the only assumption $w\in
A_{p/p_0}$ and the conclusion is  $\W_w=(p_0\,r_w,\infty)$.

\begin{remark}\rm
Observe that if $1\le p_1\le p_0\le q_0\le q_1\le \infty$ then
$$
\W_w(p_0,q_0)\subset \W_w(p_1,q_1)\subset
\W_w(1,\infty)=(r_w,\infty)=\{1<p<\infty:w\in A_p\}.
$$
\end{remark}

\begin{remark}\rm
The set $\W_w(p_0,q_0)$ can be empty: indeed, for every $1\le
p_0<q_0<\infty$, one can find $w\in A_\infty$ such that
$\W_w(p_0,q_0)=\emptyset$. A very simple example in $\re$ consists in
taking $w(x)=|x|^\alpha$ for $\alpha=q_0/p_0-1$. Note that $w\in
A_p$, $p>1$, if and only if $\alpha<p-1$ that is $p>\alpha+1$ and so
$r_w=\alpha+1$. On the other hand, $w\in RH_\infty$ and so
$s_w=\infty$. Therefore, $\W_w(p_0,q_0)=(p_0\,(1+\alpha),
q_0)=(q_0,q_0)=\emptyset$.
\end{remark}

\begin{proof}[Proof of Lemma \ref{lemma:interval-w}]
We do the case $q_0<\infty$, leaving the other one to the reader. If
$p>p_0\,r_w$ then $p/p_0>r_w$ and so $w\in A_{p/p_0}$. If,
additionally, $p<q_0/(s_w)'$ then $(q_0/p)'<s_w$ and so $w\in
RH_{(q_0/p)'}$. Therefore we have shown that $\big(\,p_0\,r_w,
q_0/(s_w)'\big)\subset \W_w(p_0,q_0)$.

To prove the converse, we observe that, by $(iii)$ in Proposition
\ref{prop:weights}, if $w\in A_{r_w}$ then $r_w=1$: if $w\in A_{r_w}$
for $r_w>1$, we have $w\in A_r$ for some $1<r<r_w$ which contradicts
the definition of $r_w$. In the same way, but this time by $(iv)$ in
Proposition \ref{prop:weights}, if $w\in RH_{s_w}$ then $s_w=\infty$.

Let $p\in \W_w(p_0,q_0)$. Since $w\in A_{p/p_0}$ then $r_w\le p/p_0$.
Besides, $r_w\neq p/p_0$ since $p/p_0>1$ and so $p>p_0\,r_{w}$. On
the other hand, $w\in RH_{(q_0/p)'}$ yields that $s_w\ge (q_0/p)'$.
Besides, $s_w\neq (q_0/p)'$ since $q_0/p>1$. This gives
$p<q_0/(s_w)'$ as desired.
\end{proof}

The duality for these classes goes as follows:

\begin{lemma}\label{lemma:duality}
Given $p_0<p<q_0$, we have
$$
w\in A_{\frac{p}{p_0}}\cap RH_{\left(\frac{q_0}{p}\right)'}
\qquad
\Longleftrightarrow
\qquad
w^{1-p'}\in A_{\frac{p'}{(q_0)'}}\cap
RH_{\left(\frac{(p_0)'}{p'}\right)'}.
$$
In other words, $p\in\W_{w}(p_0,q_0)$ if and only if $p'\in
\W_{w^{1-p'}}\big(\,(q_0)',(p_0)'\big)$.
\end{lemma}

\begin{proof}
Set $q=\big(\frac{q_0}{p}\big)'\,(\frac{p}{p_0}-1)+1$. Using
$(vi)$ and $(vii)$ in Proposition \ref{prop:weights} we have
$$
w\in A_{\frac{p}{p_0}}\cap RH_{\left(\frac{q_0}{p}\right)'}
\quad
\Longleftrightarrow
\quad
w^{\left(\frac{q_0}{p}\right)'}\in
A_{\left(\frac{q_0}{p}\right)'\,(\frac{p}{p_0}-1)+1}
=
A_q
\quad
\Longleftrightarrow
\quad
w^{\left(\frac{q_0}{p}\right)'(1-q')}\in A_{q'}
$$
and
$$
w^{1-p'}\in A_{\frac{p'}{(q_0)'}}\cap
RH_{\left(\frac{(p_0)'}{p'}\right)'}
\quad
\Longleftrightarrow
\quad
w^{(1-p')\,\left(\frac{(p_0)'}{p'}\right)'}\in
A_{\left(\frac{(p_0)'}{p'}\right)'\,\big(\,\frac{p'}{(q_0)'}-1\big)+1}.
$$
Direct computations show
$$
\left(\frac{q_0}{p}\right)'(1-q')
=
(1-p')\,\left(\frac{(p_0)'}{p'}\right)'
\qquad\mbox{and}\qquad
q'=\left(\frac{(p_0)'}{p'}\right)'\,\Big(\frac{p'}{(q_0)'}-1\Big)+1.
$$
\end{proof}

\begin{remark}\rm
Fix $1<p<\infty$. Observe that if $w$ is any given weight so that $w$, $w^{1-p'}\in
L^1_{\rm loc}$, then a given linear operator $T$ is bounded on
$L^p(w)$ if and only if its adjoint (with respect to $dx$) $T^*$ is
bounded on $L^{p'}(w^{1-p'})$. Therefore,
$$
T:L^p(w)\longrightarrow L^p(w),
\quad\mbox{for all }
w\in A_{\frac{p}{p_0}}\cap RH_{\left(\frac{q_0}{p}\right)'}
$$
if and only if
$$
T^*:L^{p'}(w)\longrightarrow L^{p'}(w),
\quad\mbox{for all }
w\in A_{\frac{p'}{(q_0)'}}\cap RH_{\left(\frac{(p_0)'}{p'}\right)'}.
$$
\end{remark}

We finish this section by giving families of weights on which $r_w$
and $s_w$ can be easily computed.

\begin{lemma}
Let $f, g\in L^1(\re^n)$ be nontrivial functions, $r\ge 1$ and
$1<s\le\infty$. Then:
\begin{list}{$(\theenumi)$}{\usecounter{enumi}\leftmargin=.8cm
\labelwidth=0.7cm\itemsep=0.2cm\topsep=.1cm
\renewcommand{\theenumi}{\roman{enumi}}
}
\item Let $w=(M f)^{-(r-1)}$ then $r_{w}=r$ and $s_{w}=\infty$, that is,
$w\in A_p\cap RH_{\infty}$ for all $p>r$ \textup{(}and $p=r$ if
$r=1$\textup{)}.

\item Let $w=(M f)^{1/s}$ then $r_{w}=1$ and $s_{w}=s$, that is,
$w\in A_1\cap RH_{q}$ for all $q<s$ \textup{(}and $q=s$ if
$s=\infty$\textup{)}.

\item If $w=(M f)^{-(r-1)}+ (M g)^{1/s}$ then
$w\in A_p\cap RH_{q}$ for all $p>r$ and $q<s$ \textup{(}and $p=r$ if
$r=1$ and $q=s$ if $s=\infty$\textup{)}. Thus,  $r_w\le r$ and
$s_w\ge s$.

\end{list}
\end{lemma}

\begin{proof}
The cases $r=1$ or $s=\infty$ are trivial. Given a nontrivial
function $f\in L^1(\re^n)$ and $\alpha\in \re$ we write
$v_\alpha=(Mf)^\alpha$. If $\alpha=0$ then $v_\alpha=1\in A_1\cap
RH_\infty$. If $0<\alpha<1$ then $v_\alpha\in A_1$ (see for instance
\cite{GR}). If $\alpha<0$ then we see that $v_\alpha\in RH_\infty$:
for a.e. $x\in B$
$$
Mf(x)^\alpha
=
\big(Mf(x)^{\frac12})^{2\,\alpha}
\lesssim
\Big(\aver{B} (Mf)^{1/2}\Big)^{2\,\alpha}
\le
\aver{B} (Mf)^{\alpha},
$$
where we have used that $(M f)^{1/2}\in A_1$ and also Jensen's
inequality for the convex function $t\mapsto t^{2\,\alpha}$. Finally,
it is easy to show that $v_\alpha\notin A_\infty$ for $\alpha\ge 1$.
Indeed, assume that $v_\alpha= (Mf)^\alpha\in A_p$ for some $1\le
p<\infty$. Then, $v_1=v_\alpha^{1/\alpha}=M f\in A_p$ as $\alpha\ge
1$. By $(vi)$ in Proposition \ref{prop:weights} we have that
$v_1^{1-p'}\in A_{p'}$ and thus $M$ is bounded on
$L^{p'}(v_1^{1-p'})$. Applying this estimate to $f\in
L^{p'}(v_1^{1-p'})$ (as $f\in L^1(\re^n)$) we obtain  that $M f\in
L^{1}(\re^n)$ which only happens when $f\equiv 0$. This leads us to a
contradiction since we have assumed that $f$ is nontrivial.

We turn to showing  $(i)$. As  $w=v_{-(r-1)}$, then  $w\in
RH_\infty$. Next, given $p>r$ the number $\alpha=(r-1)/(p-1)$
satisfies $0<\alpha<1$ and thus $v_\alpha\in A_1$. Notice that
$w=1\cdot v_\alpha^{1-p}\in A_p$ (here we are using the ``easy'' part
of the factorization of weights: if $w_1, w_2\in A_1$ then
$w_1\,w_2^{1-p}\in A_p$). This shows that $w\in A_p$ for all $p>r$
and then $r_w\le r$. To conclude we observe that $r_w=r$ as $w\notin
A_r$: otherwise we would have $w^{1-r'}= Mf\in A_{r'}$  which cannot
be the case as  seen above.

We now consider $(ii)$. Notice that $w=v_{1/s}$ with $1<s<\infty$ and thus
$w\in A_1$. Given $1<q<s$, we see that $w\in RH_q$. Note that
$w^q=v_{q/s}\in A_1$ as $q/s<1$. Then, by $(vii)$ in Proposition
\ref{prop:weights} it follows that $w\in RH_q\cap A_1$. Next, $w\notin RH_{s}$. If it were, then $w\in RH_{s+\ep}$ for some $\ep>0$ and in particular $w^s=Mf \in A_{\infty}$ which is not true.
Hence, $s_w=s$.

Note that $(iii)$ follows from $(i)$ and $(ii)$ as $w=w_1+w_2$ where
$w_1=(Mf)^{-(r-1)}\in A_p\cap RH_\infty$ and $w_2=(M f)^{1/s}\in
A_1\cap RH_q$ and $p>r$, $s<s_w$.
\end{proof}

\begin{remark}\rm
There are examples of functions $f$, $g$ for which in $(iii)$ we have
$r_w<r$ and/or $s_w> s$. For instance, if $f=g=\bigchi_{B_0}$ with
$B_0=B(0,1)$ then we have $M f(x)\approx (1+|x|)^{-n}$ and thus
$$
w(x)\approx (1+|x|)^{n\,(r-1)}\approx M f(x)^{-(r-1)}.
$$
Then, $r_w=r$ and $s_w=\infty$ (no matter the value of $s$). Similar
examples can be given in the other direction.
\end{remark}

\begin{remark}\rm
The limit case in the latter result consists of taking $f$ a Dirac
mass at some given point $x_0$, say $x_0=0$ for simplicity. In this
case $Mf(x)=c|x|^{-n}$  is a power weight. In $(i)$, $(ii)$ and
$(iii)$ we respectively have $w_1(x)=c|x|^{n\,(r-1)}$,
$w_2(x)=c|x|^{-n/s}$. Notice that $w_1\notin A_{r}$, as $w_1^{1-r'}\notin
L^1_{\rm loc}(\re^n)$. Also, $w^s\notin RH_s$ as $w^s\notin L^1_{\rm
loc}(\re^n)$.
\end{remark}

\subsection{Extrapolation}\label{subsection:extrapo}

Rubio de Francia's extrapolation theorem is a very powerful tool in
Harmonic Analysis, see \cite{Rub} and \cite{Gar}: if some given
operator $T$ is bounded on $L^{p_0}(w)$ for every $w\in A_{p_0}$ and
some $1\le p_0<\infty$, then it is bounded on $L^{p}(w)$ for all
$1<p<\infty$ and all $w\in A_p$. So, the weighted norm inequality for
one single exponent propagates to the whole range $(1,\infty)$.
Notice that in our case the natural range of exponents is no longer
$(1,\infty)$ but $(p_0,q_0)\subset (1,\infty)$.

Here we extend Rubio de Francia's result, showing that there is an
extrapolation theorem adapted to the interval $(p_0,q_0)$ which
involves the classes of weights $A_{\frac{p}{p_0}}\cap
RH_{\left(\frac{q_0}{p}\right)'}$. To state such result we first make
some reductions. As it was observed in \cite{CMP} (see also
\cite{CGMP}), one does not need to work with specific operator(s)
since nothing about the operators themselves is used (like linearity
or sublinearity) and they play no role. In other words, extrapolation
is something about weights and pairs of functions. This point of view
is very useful, for instance, when one tries to prove vector-valued
inequalities since, as we see below, they follow at once from the
corresponding scalar estimates.

So, sticking to the notation in \cite{CMP},  $\F$ denotes a family of
ordered pairs of non-negative, measurable functions $(f,g)$. In what
follows, anytime we state an estimate
$$
\|f\|_{L^p(w)}
\le
C\|g\|_{L^p(w)},
\qquad
(f,g)\in \F,
$$
we mean that it holds for all $(f,g)\in \F$ for which the left-hand
side is finite. The same is assumed when $L^{p,\infty}$ is written in
place of $L^p$ in the left hand side.

We can state our extrapolation result.

\begin{theor}\label{theor:extrapol}
Let $0<p_0<q_0\le \infty$. Suppose that there exists $p$ with $p_0\le
p\le q_0$, and $p<\infty$ if $q_0=\infty$, such that for $(f,g)\in
\F$,
\begin{equation}\label{extrapol-p}
\|f\|_{L^p(w)}
\le
C\|g\|_{L^p(w)},
\qquad
\mbox{for all }w\in A_{\frac{p}{p_0}}\cap
RH_{\left(\frac{q_0}{p}\right)'}.
\end{equation}
Then, for all $p_0<q<q_0$ and $(f,g)\in\F$ we have
\begin{equation}\label{extrapol-q}
\|f\|_{L^q(w)} \le C\,\|g\|_{L^q(w)},
\qquad
\mbox{for all }w\in A_{\frac{q}{p_0}}\cap
RH_{\left(\frac{q_0}{q}\right)'}.
\end{equation}
Moreover, for all $p_0<q,r<q_0$ and $\{(f_j,g_j)\}\subset\F$ we have
\begin{equation}\label{extrapol:v-v}
\Big\|
\Big(\sum_j (f_j)^r\Big)^{1/r}
\Big\|_{L^q(w)}
\le
C\,\Big\|
\Big(\sum_j (g_j)^r\Big)^{1/r}
\Big\|_{L^q(w)},
\quad
\mbox{for all }w\in A_{\frac{q}{p_0}}\cap
RH_{\left(\frac{q_0}{q}\right)'}.
\end{equation}
\end{theor}

The proof of this result is in Section
\ref{subsec:proof:theor:extrapol}. As an immediate consequence we can
also extrapolate from weak-type estimates:

\begin{corol}\label{corol:extrapol-weak}
Let $0<p_0<q_0\le \infty$. Suppose that there exists $p$ with $p_0\le
p\le q_0$, and $p<\infty$ if $q_0=\infty$, such that for $(f,g)\in
\F$,
\begin{equation}\label{extrapol-p:weak}
\|f\|_{L^{p,\infty}(w)}
\le
C\,\|g\|_{L^{p}(w)}
\quad
\mbox{for all }w\in A_{\frac{p}{p_0}}\cap
RH_{\left(\frac{q_0}{p}\right)'}.
\end{equation}
Then, for all $p_0<q<q_0$ and $(f,g)\in\F$ we have
\begin{equation}\label{extrapol-q:weak}
\|f\|_{L^{q,\infty}(w)}
\le
C\,\|g\|_{L^{q}(w)}
\quad
\mbox{for all }w\in A_{\frac{q}{p_0}}\cap
RH_{\left(\frac{q_0}{q}\right)'}.
\end{equation}
\end{corol}

\begin{proof}
We follow the
simple method used in \cite{GM}, for which the point of view of
 pairs of functions is particularly useful. Given $(f,g)\in \F$
and  any $\lambda>0$ we define a new pair of functions $(f_\lambda,g)$
where $f_\lambda=\lambda\,\bigchi_{E_\lambda(f)}$   and
$E_\lambda(f)=\{f>\lambda\}$. Thus \eqref{extrapol-p:weak}
implies
$$
\|f_\lambda\|_{L^p(w)}
=
\lambda\,w(E_\lambda(f))^\frac1p
\le
\sup_\lambda \lambda\,w(E_\lambda(f))^\frac1p
=
\|f\|_{L^{p,\infty}(w)}
\le
C\,\|g\|_{L^{p}(w)}
$$
for all $w\in A_{\frac{p}{p_0}}\cap RH_{\left(\frac{q_0}{p}\right)'}
$.  Applying Theorem \ref{theor:extrapol}, the family
$\widetilde{\F}$ of  pairs $(f_\lambda,g)$  satisfy \eqref{extrapol-q} with
$C$ independent of $\lambda$, and taking the supremum on
$\lambda>0$ we obtain \eqref{extrapol-q:weak}.
\end{proof}

\begin{remark}\rm
Define the following sets, given an operator $T$ defined at
least on $C_{0}^\infty(\RR^n)$:
$$
\W(T)=\big\{(p,w)\in (1,\infty)\times A_{\infty} : \|Tf\|_{L^p(w)}
\lesssim \|f\|_{L^p(w)}\big\};
$$
for $1<p<\infty$, $\W^{p}(T)=\{w\in  A_{\infty} : (p,w) \in \W(T)\}$;
and for $w \in  A_{\infty}$, $\W_{w}(T)= \{p \in (1,\infty) : (p,w)
\in \W(T)\}$.

Next define for $1\le p_{0}< q_{0}\le \infty$
$$
\W\big(\, p_0,q_0\big)
=
\big\{(p,w)\in (p_{0},q_{0})\times A_{\infty} : w\in A_{\frac{p}{p_0}}\cap
RH_{\left(\frac{q_0}{p}\right)'}
\big\};
$$
for $p_{0}<p<q_{0}$, $\W^p\big(\, p_0,q_0\big)= \{w\in A_{\infty}:
(p,w) \in \W\big(\, p_0,q_0\big) \}$ and for $w\in A_{\infty}$,
$\W_{w}\big(\, p_0,q_0\big)= \{p \in (p_{0},q_{0} ) : (p,w) \in
\W\big(\, p_0,q_0\big) \}$. Recall that the smallest $p_{0}$ (resp.
the largest $q_{0}$), the largest the class $\W\big(\, p_0,q_0\big)
$.

For example, if $T$ is a Calder\'on-Zygmund operator, then $\W(T)$
contains the largest of all classes, namely $\W(1,\infty)$ and this
is optimal. Theorem \ref{theor:main-w} (with $S =I$ and
$\D=L^{p_{0}}$) provides us with a sufficient condition on $T$ to
obtain that $\W\big(\, p_0,q_0\big)
 \subset \W(T)$.

Our extrapolation result shows that, given $T$ and $p$, if some
$\W^p(p_{0},q_{0})$ is contained in $\W^p(T)$ then  for all $q \in
(p_{0}, q_{0})$, $\W^q(p_{0}, q_{0})$ is contained in $\W^q(T)$.
In other words, $\W^p(p_{0},q_{0}) \subset \W^p(T)$ for one $p$
implies $\W(p_{0},q_{0}) \subset \W(T)$. The class of weights
$\W^p(p_{0},q_{0})$ is thus the natural one for weighted  $L^p$
boundedness within the range $p_{0}<p<q_{0}$. However, the
inclusion could be strict for a particular operator $T$ as we
will see in  \cite{AM4}.
\end{remark}

\section{Extension to spaces of homogeneous type}\label{section:SHT}

In
\cite{AM3},  we apply our results in $\re^n$ equipped with the doubling measure
$d\mu(x)=w(x)\,dx$ with $w\in A_\infty$ (in this case
$w(\re^n)=\infty$).  In \cite{AM4}, we change $\RR^n$ to a manifold or a Lie group. Hence, one needs to discuss the extension of our results to spaces of homogeneous type.

Let $(\X,d,\mu)$ be a space of homogeneous type,  that is, a set $\X$
endowed with a distance $d$ (and even a quasi-distance) and a non-negative Borel measure $\mu$ on
$\X$ such that the doubling condition
\begin{equation}\label{doubling}
\mu(B(x,2\,r)) \le C_0\, \mu(B(x,r))<\infty,
\end{equation}
holds for all $x\in\X$ and $r>0$, where $B(x,r)=\{y\in\X:d(x,y)<r\}$.

The results from Harmonic Analysis that we have used in  Euclidean
spaces remain true in this context (see for example
\cite{CW}, \cite{Chr}, \cite{Ste}). For instance, Vitali's covering
lemma,  weak-type $(1,1)$ hence strong-type $(p,p)$ for $1<p\le
\infty$ of the Hardy-Littlewood maximal function, Whitney's
covering lemma \dots The theory of Muckenhoupt weights runs
parallel to the classical case and one may prove all the
statements in Proposition \ref{prop:weights} with the appropriate
changes (see \cite[Chapter I]{ST}).

Hence, Theorems \ref{theor:good-lambda:w},  \ref{theor:main-w},
\ref{theor:ACDH:w}, \ref{theor:shen}, \ref{theor:comm},
\ref{theor:extrapol} all have their counterpart in spaces of
homogeneous type with  almost identical proofs whenever
$\mu(\X)=\infty$.

When $\mu(\X)<\infty$ (for example, $\X$ is a bounded Lipschitz
domain in $\RR^n$) some adjustments are needed. In Theorem
\ref{theor:good-lambda:w}, assuming that $F\in L^1$ then the two
parameter good-$\lambda$ estimate \eqref{good-lambda:w} holds for
$\lambda>\lambda_0=C_0\,\mu(\X)^{-1}\,(\|F\|_{L^1}+\|H_1\|_{L^1})$.
This condition guarantees that $\mu(E_\lambda)<\mu(\X)$ and so
$E_\lambda\subsetneq \X$. The Whitney covering argument can be
performed and the proof presented above works in the same way. Thus,
when proving the analog of \eqref{good-lambda:Lp:w}, one has to split
the integral in two parts: $\lambda\ge \lambda_0$ and $\lambda\le
\lambda_0$. For the first one, we use \eqref{good-lambda:w}. The
piece  $\lambda\le \lambda_0$ is estimated by observing that $w\{M
F>\lambda\}\le w(\X)<\infty$ (since $\mu(\X)<\infty$ if and only if
$\X$ is bounded, see for instance \cite{Ma2}). Thus, it can be proved
that
$$
\|M F\|_{L^p(w)}
\le
C\,\big(\|G\|_{L^p(w)} + \|M
H_1\|_{L^p(w)}+\|H_2\|_{L^p(w)}+\|F\|_{L^1(\mu)}+\|MH_1\|_{L^1(\mu)}\big).
$$
The same occurs with the estimates in $L^{p,\infty}(w)$.

The latter inequality allows one to obtain Theorem \ref{theor:main-w}
assuming further that $T$ is bounded on $L^{p_0}$ (this happens all
the time in applications, see \cite{AM3}). The only change is for
 the term $\|F\|_{L^1}$ where $F=\vert Tf\vert ^{p_{0}}$ (notice that $H_1=H_{2}=0$ in this case):
$$
\|F\|_{L^1(\mu)}
=
\|T f\|_{L^{p_0}(\mu)}^{p_0}
\lesssim
\|f\|_{L^{p_0}(\mu)}^{p_0}
\le
\|f\|_{L^{p}(w)}^{p_0}\,\int_\X w^{1-(p/p_0)'} d\mu
\lesssim
\|f\|_{L^{p}(w)}^{p_0}.
$$
For the last inequality, we observe that since  $w\in A_{(p/p_0)}$,
then $w^{1-(p/p_0)'}\in A_{(p/p_0)'}$ and so it is a doubling measure
which implies as noted before that $w^{1-(p/p_0)'}(\X)<\infty$ as
$\X$ is bounded. Similar modifications can be carried out with
Theorems \ref{theor:ACDH:w} and \ref{theor:shen}. Precise statements
and details of  proofs are left to the interested reader.

\section{Proofs of the main results}\label{section:proofs1}

 We prove Theorem \ref{theor:good-lambda:w}, Theorem \ref{theor:comm}, Theorem \ref{theor:extrapol}.

\subsection{Proof of Theorem \ref{theor:good-lambda:w}}
\label{subs:proof:theor:good-lambda:w}

The proof follows the ideas in \cite{Aus-mem}. It suffices to
consider the case $H_2=G$: indeed, set $\widetilde{G}=G+H_2$. Then
\eqref{H-Q} holds with $\widetilde{G}$ in place of $H_2$ and also
\eqref{G-Q} holds with $\widetilde{G}$ in place of $G$.

So from now on we assume that $H_2=G$. Set  $E_\lambda=\{MF+ M
H_1>\lambda\}$ which is assumed to have finite measure (otherwise
there is nothing to prove). As $M$ is the uncentered maximal function
(over cubes instead of balls),  $E_\lambda$ is an open set. Hence,
Whitney's decomposition gives us a family of pairwise disjoint cubes
$\{Q_j\}_j$ so that $E_\lambda=\cup_j Q_j$ and with the property that
$4\,Q_j$ meets $E_\lambda^c$, that is, there exists $x_j\in 4\,
Q_{j}$ such that
$$
M F(x_j)+ M H_1(x_j)\le \lambda.
$$
Set $B_\lambda=\{MF>K\,\lambda, 2\,G\le \gamma\,\lambda\}$. Since
$K\ge 1$ we have that $B_\lambda\subset E_\lambda$. Therefore
$B_\lambda\subset \cup_j B_\lambda\cap Q_j$. For each $j$ we assume
that $B_\lambda\cap Q_j\neq \emptyset$ (otherwise we discard this
cube) and so there is $\bar{x}_j\in Q_j$ so that $G(\bar{x}_j)\le
\gamma\,\lambda/2$. Since $M F(x_j)\le \lambda$, there is $C_{0}$
depending only on dimension such that  for every $K\ge C_0$ we have
\begin{align*}
\lefteqn{\hskip-0.5cm
|B_\lambda\cap Q_j|
\le
\big| \{M F>K\, \lambda\} \cap Q_{j}\big|
\le
\big|\{M (F\bigchi_{8\,Q_j})>(K/C_0)\, \lambda\}\big|
}
\\
&\le
\big|\{M (G_{8\,Q_j}\bigchi_{8\,Q_j})>(K/2\,C_0)\, \lambda\}\big| +
\big|\{M (H_{8\,Q_j}\bigchi_{8\,Q_j})>(K/2\,C_0)\, \lambda\}\big|,
\end{align*}
where we have used $F\bigchi_{8\,Q_j}\le
G_{8\,Q_j}\bigchi_{8\,Q_j}+H_{8\,Q_j}\bigchi_{8\,Q_j}$ a.e.   and
$\bigchi_{8\,Q_j}$ is the indicator function of $8\, Q_{j}$. Let
$c_{p}$ be the weak-type $(p,p)$ bound of the maximal function. By
\eqref{G-Q} and $\bar{x}_j\in Q_j\subset 8\,Q_j$, we obtain
\begin{align*}
\big|\{M (G_{8\,Q_j}\bigchi_{8\,Q_j})>(K/2\,C_0)\,
\lambda\}\big|
&\le
\frac{2\,C_0\,c_1}{K\,\lambda}\, \int_{8\,Q_j} G_{8\,Q_j}
\le
\frac{2\,C_0\,c_1}{K\,\lambda} \,|8\,Q_j|\, G(\bar{x}_j)
\\
&\le
\frac{8^{n}\,C_0\,c_1}{K}\, |Q_j| \,\gamma.
\end{align*}
Next, assume first that $q<\infty$. By  \eqref{H-Q} and  $x_j$,
$\bar{x}_j\in 8\,Q_j$, we obtain
\begin{align*}
\lefteqn{\hskip-.5cm
\big|\{M (H_{8\,Q_j}\bigchi_{8\,Q_j})>(K/2\,C_0)\,
\lambda\}\big|
\le
\left(\frac{2\,C_0\,c_q}{K\,\lambda}\right)^q\, \int_{8\,Q_j}
H_{8\,Q_j}^q}
\\
&
\le
\left(\frac{2\,C_0\,c_q}{K\,\lambda}\right)^q\,|8\,Q_j|\,a^q\,\big(
M F(x_j) + M H_1(x_j)+G(\bar{x}_j)\big)^q
\le
\frac{(4\,C_0\,c_q\,a)^q\,8^{n}}{K^q}\, |Q_j|.
\end{align*}
These two estimates yield
$$
|B_\lambda\cap Q_j|
\le
C\,\left(\frac{a^q}{K^q}+\frac{\gamma}{K}\right)\,|Q_j|.
$$
At this point, we use that $w\in RH_{s'}$. If $s'<\infty$, for any
cube $Q$ and any measurable set $E\subset Q$ we have
$$
\frac{w(E)}{w(Q)}
\le
\frac{|Q|}{w(Q)}\,
\Big(\aver{Q}
w^{s'}\Big)^{\frac1{s'}} \,\left(\frac{|E|}{|Q|}\right)^{\frac1s}
\le
C_w\,\left(\frac{|E|}{|Q|}\right)^{\frac1s}.
$$
Note that the same conclusion holds in the case $s'=\infty$. Applying
this to $B_\lambda\cap Q_j\subset Q_j$ we have
$$
w(B_\lambda\cap Q_j)
\le
C_w\,C\,\left(\frac{a^q}{K^q}+\frac{\gamma}{K}\right)^{\frac1s}\,w(Q_j).
$$
Hence, using that the Whitney cubes are disjoint we have
$$
w(B_\lambda)
\le
\sum_j w(B_\lambda\cap Q_j)
\le
C\,\left(\frac{a^q}{K^q}+\frac{\gamma}{K}\right)^{\frac1s}\,\sum_j
w(Q_j)
=
C\,\left(\frac{a^q}{K^q}+\frac{\gamma}{K}\right)^{\frac1s}\,w(E_\lambda)
$$
which is \eqref{good-lambda:w}.

When $q=\infty$, then by \eqref{H-Q}
$$
\|M (H_{8\,Q_j}\bigchi_{8\,Q_j})\|_{L^\infty}
\le
\|H_{8\,Q_j}\bigchi_{8\,Q_j})\|_{L^\infty}
\le
a\,\big(M F(x_j) + MH_1(x_j)+G(\bar{x}_j)\big)
\le 2\,a\, \lambda.
$$
Thus choosing $K\ge 4\,a\,C_0$ it follows that $\{M
(H_{8\,Q_j}\bigchi_{8\,Q_j})>(K/2\,C_0)\, \lambda\}=\emptyset$.
Proceeding as before, we get the desired estimate (with $K^{-q}=0$).

Next we show \eqref{good-lambda:Lp:w} when it is assumed that $MF\in
L^p(w)$. Integrating the two-parameter good-$\lambda$ inequality
\eqref{good-lambda:w} against $p\,\lambda^{p-1}\, d\lambda$ on
$(0,\infty)$, for $0<p<\infty$,
$$
\|M F\|_{L^p(w)}^p
\le
C\,K^p\,
\left(
\frac{a^q}{K^q}+\frac{\gamma}{K}\right)^{\frac1s} \,\big(
\|MF\|_{L^p(w)}^p +\|M H_1\|_{L^p(w)}^p\big)+
\frac{2^p\,K^p}{\gamma^p}\,\|G\|_{L^p(w)}^p.
$$
Thus, as $\|M F\|_{L^p(w)}<\infty$, for $0<p<\frac{q}{s}$ we can
choose $K$ large enough and then $\gamma$ small enough so that the
constant in front of the first term in the right-hand side is smaller
than $\frac12$,  leading us to \eqref{good-lambda:Lp:w}. In the same
way, but this time assuming that $MF\in L^{p,\infty}(w)$, one shows
the corresponding estimate in $L^{p,\infty}(w)$ .

Observe that in the case $q=\infty$, $K$ is already chosen and we
only have to take some small $\gamma$. Thus, the corresponding
estimates holds for $0<p<\infty$ no matter the value of $s$.

Now, we consider the case $p\ge 1$ and $F\in L^1$. We assume that
the right-hand side of \eqref{good-lambda:Lp:w} is finite,
otherwise there is nothing to prove. It suffices to consider the
case $w\in L^\infty$: indeed we can take $w_N=\min\{w,N\}$ with
$N>0$. As $w\in RH_{s'}$ then $w_N\in RH_{s'}$ with constant that
is uniformly controlled in $N$. Notice that if we show
\eqref{good-lambda:Lp:w} with $w_N$ and with constants that do not
depend on $N$, by taking limits as $N\rightarrow \infty$, we
conclude the desired estimate with $w$.

So we assume that $w\in L^\infty$. Let $f$  be the
non-negative function defined by
$f(\lambda)=p\lambda^{p}\,w\{MF>\lambda\}$, $\lambda>0$. Notice that for any $0<\lambda_0<\lambda_1<\infty$, $\int_{\lambda_0}^{\lambda_1}
f(\lambda)\,\frac{d\lambda}{\lambda}$ exists and is finite.  By
\eqref{good-lambda:w}  we have
\begin{eqnarray*}
\lefteqn{\int_{\lambda_0}^{\lambda_1}
f(\lambda)\,\frac{d\lambda}{\lambda}
=
\int_{\frac{\lambda_0}{K}}^{\frac{\lambda_1}{K}}
f(K\,\lambda)\,\frac{d\lambda}{\lambda}}
\\
&\le&
C\,K^p\, 2^p
\left(
\frac{a^q}{K^q}+\frac{\gamma}{K}\right)^{\frac1s}\,\Big(
\int_{\frac{\lambda_0}{2\,K}}^{\frac{\lambda_1}{2\,K}}
f(\lambda)\,\frac{d\lambda}{\lambda} + \|M H_1\|_{L^p(w)}^p\Big)+
\frac{2^p\,K^p}{\gamma^p}\,\|G\|_{L^p(w)}^p
\\
&\le&
\frac12\,\int_{\frac{\lambda_0}{2\,K}}^{\frac{\lambda_1}{K}}
f(\lambda)\,\frac{d\lambda}{\lambda} + R
\end{eqnarray*}
where in the last inequality we have picked $K$ large enough and then
$\gamma$ small enough so that the constant in front of the first term
in the right-hand side is smaller than $1/2$. Also we have written
$R$ for the remainder terms, that is, $R=C\,\big(\|M H_1\|_{L^p(w)}^p
+\|G\|_{L^p(w)}^p\big)<\infty$. We take $\lambda_0=K^{-n}$ and
$\lambda_1=K^{m}$ with $n,m\ge 1$ and so
\begin{eqnarray*}
\int_{K^{-n}}^{K^{m-1}} f(\lambda)\,\frac{d\lambda}{\lambda}
&\le&
\int_{K^{-n}}^{K^m} f(\lambda)\,\frac{d\lambda}{\lambda}
\le
\frac12\,\int_{\frac{K^{-n-1}}{2}}^{K^{m-1}}
f(\lambda)\,\frac{d\lambda}{\lambda} + R
\\
&\le&
\frac12\,\int_{K^{-n}}^{K^{m-1}} f(\lambda)\,\frac{d\lambda}{\lambda}
+ \frac12\, \int_{\frac{K^{-n-1}}{2}}^{K^{-n}}
f(\lambda)\,\frac{d\lambda}{\lambda}+ R.
\end{eqnarray*}
Hence,
$$
\int_{K^{-n}}^{K^{m-1}} f(\lambda)\,\frac{d\lambda}{\lambda}
\le
\int_{\frac{K^{-n-1}}{2}}^{K^{-n}}
f(\lambda)\,\frac{d\lambda}{\lambda}+ 2\,R.
$$
Since $M$ is of weak-type $(1,1)$,  $w\in L^\infty$ and
 $K\ge 1$ we have
$$
\int_{\frac{K^{-n-1}}{2}}^{K^{-n}}
f(\lambda)\,\frac{d\lambda}{\lambda}
\le
C\,\|w\|_{L^\infty}\,\|F\|_{L^1}\,
\left\{
\begin{array}{l@{\hskip.5cm}l}
\log 2\,K & \mbox{ if } p=1,
\\[0.3cm]
1 & \mbox{ if } p>1,
\end{array} \right.
$$
bound which does not depend on $n$. We conclude that
$$
\|M F\|_{L^p(w)}^p
=
\int_0^\infty f(\lambda)\,\frac{d\lambda}{\lambda}
=
\lim_{n,m\rightarrow\infty} \int_{K^{-n}}^{K^{m-1}}
f(\lambda)\,\frac{d\lambda}{\lambda}
<\infty,
$$
so that  $MF\in L^p(w)$. Therefore,  \eqref{good-lambda:Lp:w} holds
with constants that do {\em not} depend on $\|w\|_{L^\infty}$. A very
similar argument applies for the weak-type estimate. Details are left
to the reader. \qed

\subsection{Proof of Theorem \ref{theor:comm}}\label{subsec:proof:theor:comm}

Before starting the proof,  let us introduce some notation (see
\cite{BS} for more details). Let $\phi$ be a Young function: $\phi
:[0,\infty)\longrightarrow [0,\infty)$ is  continuous, convex,
increasing and satisfies $\phi(0+)=0$, $\phi(\infty)=\infty$. Given a
cube $Q$ we define the localized Luxemburg's norm
$$
\|f\|_{\phi ,Q}
=
\inf
\bigg\{
\lambda>0: \aver{Q} \phi
\left(\frac{|f|}{\lambda}\right)\, \leq 1
\bigg\},
$$
and then the maximal operator
$$
M_{\phi}f(x)
=
\sup_{Q\ni x} \|f\|_{\phi ,Q}.
$$
In the definition of $\|\cdot\|_{\phi,Q}$, if the probability measure
$dx/|Q|$ is replaced by $dx$  and $Q$ by $\re^n$, then one has the
Luxemburg's norm $\|\cdot\|_{\phi}$ which allows one to define the
Orlicz space $L^{\phi}$.

Some specific examples  needed here are $\phi(t)\approx e^{t^r}$ for
$t\ge 1$ which gives the classical space ${\rm exp}  L^r$ and
$\phi(t)=t\,(1+\log^+ t)^\alpha$ with $\alpha>0$ that gives the space
$L\,(\log L)^\alpha$. In this latter case, it is well known that for
$k\ge 1$, we have $M_{L(\log L)^{k-1}} f\approx M^{k} f$ where
$M^{k}$ is the $k$-iteration of $M$.

John-Nirenberg's inequality implies that for any function $b\in \BMO$
and any cube $Q$ we have $\|b-b_Q\|_{{\rm exp} L,Q}\lesssim
\|b\|_{\BMO}$. This yields the following estimates: First,  for each
cube $Q$ and $x\in Q$
\begin{align}
\lefteqn{\hskip-1.5cm
\aver{Q} |b-b_Q|^{k\,p_0}\,|f|^{p_0}
\le
\|b-b_Q\|_{{\rm exp} L,Q}^{k\,p_0}\,\big\| |f|^{p_0}
\big\|_{L\,\,(\log L)^{k\,p_0},Q}}
\nonumber
\\
&
\lesssim
\|b\|_{\BMO}^{k\,p_0}\, M_{L\,\,(\log L)^{k\,p_0}}\big( |f|^{p_0})(x)
\lesssim
\|b\|_{\BMO}^{k\,p_0}\, M^{[k\,p_0]+2}\big( |f|^{p_0})(x),
\label{BMO-M}
\end{align}
where $[s]$ is the integer part of $s$ (if $k\,p_0\in\NN$, then one
can take $M^{[k\,p_0]+1}$). Second, for each $j\ge 1$ and each $Q$,
\begin{align}
\|b-b_{2\,Q}\|_{{\rm exp} L,2^{j}\,Q}
&\le
\|b-b_{2^{j}\,Q}\|_{{\rm exp} L,2^{j}\,Q} +|b_{2^{j}\,Q}-b_{2\,Q}|
\lesssim
\|b\|_{\BMO}+ \sum_{l=1}^{j-1}|b_{2^{l+1}\,Q}-b_{2^l\,Q}|
\nonumber
\\
&\lesssim
\|b\|_{\BMO}+ \sum_{l=1}^{j-1} \aver{2^{l+1}\,Q} |b-b_{2^{l+1}\,Q}|
\lesssim j\,\|b\|_{\BMO} \label{BMO-log}.
\end{align}

The following auxiliary result allows us to assume further that $b\in
L^\infty$. The proof is postponed until the end of this section.

\begin{lemma}\label{lemma:comm-apriori}
Let $1\le p_0<p<\infty$, $k\in\NN$ and $w\in A_\infty$. Let $T$ be a
sublinear operator bounded on $L^{p_0}$.
\begin{list}{$(\theenumi)$}{\usecounter{enumi}\leftmargin=.8cm
\labelwidth=0.7cm\itemsep=0.2cm\topsep=.3cm
\renewcommand{\theenumi}{\roman{enumi}}
}
\item If $b\in \BMO\cap L^\infty$ and $f\in
L^\infty_c$, then $T_b^k f\in L^{p_0}$.

\item Assume that for any $b\in \BMO\cap L^\infty$ and for any $f\in L^\infty_c$
we have that
\begin{equation}\label{lemma:comm-est}
\|T_b^k f\|_{L^p(w)}
\le C_0\, \|b\|_{\BMO}^k\,\|f\|_{L^p(w)},
\end{equation}
where $C_0$ does not depend on $b$ and  $f$. Then for all $b\in
\BMO$, \eqref{lemma:comm-est} holds with constant $2^k\,C_0$ instead of $C_{0}$.
\end{list}
\end{lemma}

Part $(ii)$ in this latter result ensures that it suffices to
consider the case $b\in L^\infty$ (provided the constants obtained do
not depend on $b$). So from now on we assume that $b\in L^\infty$ and
obtain \eqref{lemma:comm-est} with $C_0$ independent of $b$ and $f$.
Note that by homogeneity we can also assume that $\|b\|_{\BMO}=1$.

We proceed by induction. As mentioned in Remark \ref{remark:indu-0},
the case $k=0$ follows from Theorem \ref{theor:main-w}. We write the
case $k=1$ in full detail  and indicate
how to pass from $k-1$ to $k$ as the argument  is essentially the
same. Let us fix $p_0<p<q_0$ and $w\in A_{\frac{p}{p_0}}\cap
RH_{\left(\frac{q_0}{p}\right)'}$. We assume that $q_0<\infty$, for
$q_0=\infty$ the main ideas are the same and details are left to the
interested reader.

\

\noindent \textit{Case $k=1$}: We combine the ideas in the proof of
Theorem \ref{theor:main-w} with techniques for commutators, see
\cite{Per}. Let $f\in L^\infty_c$ and set $F=|T_b^1 f|^{p_0}$. Note
that $F\in L^1$ by $(i)$ in Lemma \ref{lemma:comm-apriori} (this is
the only place in this step where we use  that $b\in L^\infty$).
Given a cube $Q$, we set $f_{Q,b}=(b_{4\,Q}-b)\,f$ and decompose
$T_b^1$ as follows:
\begin{eqnarray*}
|T_b^1 f(x)|
&=&
|T\big((b(x)-b)\,f\big)(x)|
\le
|b(x)-b_{4\,Q}|\,|T f(x)|+|T\big((b_{4\,Q}-b)\,f\big)(x)|
\\
&\le&
|b(x)-b_{4\,Q}|\,|T f(x)|+|T(I-\A_{r(Q)})f_{Q,b}(x)| +
|T\A_{r(Q)}f_{Q,b}(x)|.
\end{eqnarray*}
With the notation of Theorem \ref{theor:good-lambda:w},  we observe
that $F\le G_Q+H_Q$ where
$$
G_Q
=
4^{p_0-1}\big(G_{Q,1}+G_{Q,2}\big)
=
4^{p_0-1}\,\big(|b-b_{4\,Q}|^{p_0}\,|T
f|^{p_0}+|T(I-\A_{r(Q)})f_{Q,b}|^{p_0}\big)
$$
and $H_Q=2^{p_0-1}\,|T\A_{r(Q)}f_{Q,b}|^{p_0}$.

We first estimate the average of  $G_{Q}$ on $Q$. Fix any $x\in Q$.
By  \eqref{BMO-M} with $k=1$,
$$
\aver{Q} G_{Q,1}
=
\aver{Q} |b-b_{4\,Q}|^{p_0}\,|T f|^{p_0}
\lesssim
\|b\|_{\BMO}^{p_0}\, M^{[p_0]+2}\big(|T f|^{p_0})(x).
$$
Using \eqref{T:I-A:comm}, \eqref{BMO-M} and \eqref{BMO-log},
\begin{align*}
\Big(\aver{Q} G_{Q,2}\Big)^{\frac1{p_0}}
&=
\Big(\aver{Q} |T(I-\A_{r(Q)})f_{Q,b}|^{p_0}\Big)^{\frac1{p_0}}
\lesssim
\sum_{j=1}^\infty \alpha_j
\Big(\aver{2^{j+1}\,Q}|f_{Q,b}|^{p_0}\Big)^{\frac1{p_0}}
\\
&
\le
\sum_{j=1}^\infty \alpha_j\, \|b-b_{4\,Q}\|_{{\rm exp}
L,2^{j+1}\,Q}\, M^{[p_0]+2}\big(|f|^{p_0})^{\frac1{p_0}}(x)
\\
&\lesssim
\|b\|_{\rm BMO}
\,M^{[p_0]+2}\big(|f|^{p_0})(x)^{\frac1{p_0}}\,\sum_{j=1}^\infty
\alpha_j\,j
\lesssim
M^{[p_0]+2}\big(|f|^{p_0})^{\frac1{p_0}}(x),
\end{align*}
since $\sum_j \alpha_j\,j<\infty$. Hence, for any $x\in Q$
$$
\aver{Q} G_{Q}
\le
C\,
\big(M^{[p_0]+2}\big(|T f|^{p_0})(x) + M^{[p_0]+2}\big(|f|^{p_0})(x)
\big)
\equiv G(x) .
$$

We next estimate the average of  $H_Q^q$ on $Q$ with $q=q_0/p_0$.
Using \eqref{T:A:comm} and proceeding as before
\begin{align*}
\lefteqn{\hskip-1.2cm
\Big(\aver{Q} H_{Q}^{q}\Big)^{\frac1{q_0}}
=
2^{(p_0-1)/p_0}\Big(\aver{Q} |T \A_{r(Q)}
f_{Q,b}|^{q_0}\Big)^{\frac1{q_0}}
\lesssim
\sum_{j=1}^\infty \alpha_j \Big(\aver{2^{j+1}\,Q} |T
f_{Q,b}|^{p_0}\Big)^{\frac1{p_0}}}
\\
&\le
\sum_{j= 1}^\infty \alpha_j \Big(\aver{2^{j+1}\,Q} |T_b^1
f|^{p_0}\Big)^{\frac1{p_0}} + \sum_{j\ge 1} \alpha_j
\Big(\aver{2^{j+1}\,Q} |b-b_{4\,Q}|^{p_0} |T f|^{p_0}\Big)^{\frac1{p_0}}
\\
&\lesssim
(M F)^{\frac1{p_0}}(x)+ \sum_{j=1}^\infty \alpha_j\,
\|b-b_{4\,Q}\|_{{\rm exp} L,2^{j+1}\,Q}\, M^{[p_0]+2}\big(|T
f|^{p_0})^{\frac1{p_0}}(\bar{x})
\\
&\lesssim
(M F)^{\frac1{p_0}}(x)+ M^{[p_0]+2}\big(|T
f|^{p_0})^{\frac1{p_0}}(\bar{x}) \,\sum_{j=1}^\infty \alpha_j\,j
\\
&
\lesssim
(M F)^{\frac1{p_0}}(x)+ M^{[p_0]+2}\big(|T
f|^{p_0})^{\frac1{p_0}}(\bar{x}),
\end{align*}
for any $x$, $\bar{x}\in Q$, where we have used that $\sum_j
\alpha_j\,j<\infty$. Thus we have obtained
$$
\Big(\aver{Q} H_{Q}^{q}\Big)^{\frac1{q}}
\le
C\,\big( MF(x)+  M^{[p_0]+2}\big(|T f|^{p_0})(\bar{x})\big) \equiv
C\, \big( M F(x)+ H_2(\bar{x})\big).
$$

As mentioned before $F\in L^{1}$.  Since $w\in RH_{(q_0/p)'}$,
applying Theorem \ref{theor:good-lambda:w} and Remark
\ref{remark:up-endpoint} (since $q_0<\infty$ implies $q<\infty$) with
$p/p_0$ in place of $p$ and $s=q_0/p$, we obtain
\begin{align*}
\|T_b^1 f\|_{L^p(w)}^{p_0}
&\le
\|M F\|_{L^\frac{p}{p_0}(w)}
\lesssim
\|G\|_{L^\frac{p}{p_0}(w)}+\|H_2\|_{L^\frac{p}{p_0}(w)}
\\
&\lesssim
\big\|M^{[p_0]+2}\big(|f|^{p_0}\big)\big\|_{L^\frac{p}{p_0}(w)}
+
\big\|M^{[p_0]+2}\big(|T f|^{p_0}\big)\big\|_{L^\frac{p}{p_0}(w)}
\\
&\lesssim
\|f\|_{L^p(w)}^{p_0}+\|T f\|_{L^p(w)}^{p_0}
\lesssim
\|f\|_{L^p(w)}^{p_0} ,
\end{align*}
where we have used the boundedness of $M$ (hence, $M^2, M^3, \dots$)
on $L^{\frac{p}{p_0}}(w)$ as $w\in A_{p/p_0}$ with $p_0<p$, and also
Remark \ref{remark:indu-0}. Let us emphasize  that none of the
constants depend on $b$ or $f$.

\

\noindent \textit{Case $k$}:
 We now sketch the induction argument.
Assume that we have already proved the cases $m=0,\dots,k-1$. Let
$f\in L^\infty_c$. Given a cube $Q$, write
$f_{Q,b}=(b_{4\,Q}-b)^k\,f$ and decompose $T_b^k$ as follows:
\begin{align*}
|T_b^{k}f(x)|
&=
|T\big((b(x)-b)^k\,f\big)(x)|
\\
&\le
\sum_{m=0}^{k-1}C_{k,m}|b(x)-b_{4\,Q}|^{k-m} |T_b^{m} f(x)| +
|T\big((b_{4\,Q}-b)^k f\big)(x)|
\\
&\lesssim
\sum_{m=0}^{k-1}|b(x)-b_{4\,Q}|^{k-m} |T_b^{m} f(x)| +
|T(I-\A_{r(Q)})f_{Q,b}(x)| + |T\A_{r(Q)}f_{Q,b}(x)|.
\end{align*}
Following the notation of Theorem \ref{theor:good-lambda:w}, we set
$F=|T_b^k f|^{p_0}\in L^1$ by $(i)$ in Lemma
\ref{lemma:comm-apriori}. Observe that $F\le G_Q+H_Q$ where
\begin{eqnarray*}
G_Q = 4^{p_0-1}\,C\,
\Big(\Big(\sum_{m=0}^{k-1}|b-b_{4\,Q}|^{k-m} |T_b^{m} f|\Big)^{p_0}
+|T(I-\A_{r(Q)})f_{Q,b}|^{p_0}\Big)
\end{eqnarray*}
and $H_Q=2^{p_0-1}\,|T\A_{r(Q)}f_{Q,b}|^{p_0}$.  Proceeding as before
we obtain for any $x\in Q$
$$
\aver{Q} G_{Q}
\le
C\,\Big(\sum_{m=0}^{k-1} M^{[(k-m)\,p_0]+2}\big(|T_b^m f|^{p_0})(x) +
M^{[k\,p_0]+2}\big(|f|^{p_0}\big)(x)\Big) \equiv G(x) ,
$$
and for  $q=q_0/p_0$
$$
\Big(\aver{Q} H_{Q}^{q}\Big)^{\frac1{q}}
\le
C\,\big( M F(x)+ \sum_{m=0}^{k-1} M^{[(k-m)\,p_0]+2}\big(|T_b^m
f|^{p_0}\big)(\bar{x})\big) \equiv C\, \big( MF(x)+
H_2(\bar{x})\big).
$$
Therefore, as $F\in L^{1}$, Theorem \ref{theor:good-lambda:w} gives
us as before
\begin{align*}
\|T_b^k f\|_{L^p(w)}^{p_0}
&\le
\|M F\|_{L^\frac{p}{p_0}(w)}
\lesssim
\|G\|_{L^\frac{p}{p_0}(w)}+\|H_2\|_{L^\frac{p}{p_0}(w)}
\\
&\lesssim
\big\|M^{[k\,p_0]+2}\big(|f|^{p_0}\big)\big\|_{L^\frac{p}{p_0}(w)}
+ \sum_{m=0}^{k-1}
\big\|M^{[(k-m)\,p_0]+2}\big(|T_b^m f|^{p_0}\big)\big\|_{L^\frac{p}{p_0}(w)}
\\
&\lesssim
\|f\|_{L^p(w)}^{p_0}+\sum_{m=0}^{k-1} \|T_b^m f\|_{L^p(w)}^{p_0}
\lesssim
\|f\|_{L^p(w)}^{p_0},
\end{align*}
where we have used the boundedness on $L^\frac{p}{p_0}(w)$ of the
iterations of $M$ (as $w\in A_{p/p_0}$ and $p>p_0$) and   the
induction hypothesis on   $T_b^m$, $m=0,\dots,k-1$. Let us point out again
that none of the constants involved in the proof depend on $b$ and
$f$.

\begin{proof}[Proof of Lemma \ref{lemma:comm-apriori}]
Some of the ideas of the following argument are taken from
\cite{Per} where this is proved for Calder\'on-Zygmund operators.
Note that there, one has size and smoothness estimates for the
kernels and here such  conditions are not assumed.

Fix $f\in L^\infty_c$. Note that $(i)$ follows easily observing that
$$
|T_b^k f(x)|
\lesssim
\sum_{m=0}^k |b(x)|^{m-k}\,|T(b^m\,f)(x)|
\le
\|b\|_{L^\infty}\sum_{m=0}^k |T(b^m\,f)(x)| \in L^{p_0},
$$
since $b\in L^\infty$, $f\in L^\infty_c$ imply that $b^m\,f\in
L^\infty_c\subset L^{p_0}$ and, by assumption, $T(b^m\,f)\in
L^{p_0}$.

To obtain $(ii)$, we fix $b\in \BMO$ and $f\in L^\infty_c$. Let $Q_0$
be a cube such that $\supp f\subset Q_0$. We may assume that
$b_{Q_0}=0$ since otherwise we can work with
$\widetilde{b}=b-b_{Q_0}$ and observe that
$T_b^k=T_{\widetilde{b}}^k$ and
$\|b\|_{\BMO}=\|\widetilde{b}\|_{\BMO}$. Note that  for all $m=0,
\ldots, k$, we have that $|b^m\, f|$ and $\big|T(b^m\,f)\big|$ are
finite almost everywhere since they belong  to  $L^{p_0}$.

Let $N>0$ and define $b_N$ as follows: $b_N(x)=b(x)$ when $-N\le
b(x)\le N$, $b_N(x)=N$ when $b(x)>N$ and $b(x)=-N$ when $b(x)<-N$.
Then, it is immediate to see that $|b_N(x)-b_N(y)|\le |b(x)-b(y)|$
for all $x$, $y$.  Thus, $\|b_N\|_{\BMO}\le 2\,\|b\|_{\BMO}$. As
$b_N\in L^\infty$ we can use \eqref{lemma:comm-est} and
$$
\|T_{b_N}^k f\|_{L^p(w)}
\le
C_0\, \|b_N\|_{\BMO}^k\,\|f\|_{L^p(w)}
\le
C_0\,2^k\, \|b\|_{\BMO}^k\,\|f\|_{L^p(w)}<\infty.
$$
To conclude, by Fatou's lemma, it suffices to show that $|T_{b_{N_j}}
f(x)|\longrightarrow |T_b^k f(x)|$ for a.e. $x\in \re^n$ and for some
subsequence $\{N_j\}_j$ such that $N_j\rightarrow \infty$.

As $|b_N|\le |b|\in L^{p}(Q_0)$ for any $1\le p<\infty$,  the
dominated convergence theorem yields that $(b_N)^m\,f\longrightarrow
b^m\,f$ in $L^{p_0}$ as $N\rightarrow\infty$ for all $m=0,\dots,k$.
Therefore, as $T$ is bounded on $L^{p_0}$ it follows that
$T\big((b_N)^m\,f-b^m\,f\big)\longrightarrow 0$ in $L^{p_0}$. Thus,
there exists a subsequence $N_j\rightarrow \infty$ such that
$T\big((b_{N_j})^m\,f-b^m\,f\big)(x)\longrightarrow 0$ for a.e. $x\in
\re^n$ and for all $m=1,\dots,k$. In this way we obtain
\begin{align*}
\lefteqn{\big| |T_{b_{N_j}}^k f(x)|-|T_{b}^k f(x)|
\big|
\lesssim
\big|
T
\big(
\big[(b_{N_j}(x)-b_{N_j})^k-(b(x)-b)^k\big]\,f
\big)(x)
\big|}
\\
&\lesssim
\sum_{m=0}^k
|b_{N_j}(x)|^{k-m}\,\big|T\big((b_{N_j})^m\,f-b^m\,f\big)(x)\big| +
\big|b_{N_j}(x)^{k-m}-b(x)^{k-m}\big|\,\big|T(b^m\,f)(x)\big|
\end{align*}
and as desired we get that $|T_{b_{N_j}} f(x)|\longrightarrow |T_b^k
f(x)|$ for a.e. $x\in \re^n$.
\end{proof}

\begin{remark}\label{ref:proof::multi-comm:p-big}\rm
The proof just finished can be adapted to the situation of
multilinear commutators with no much effort. We just sketch some of
the ideas leaving the details to the interested reader. Let us
introduce some notation. Given $\vec{b}=(b_1,\dots,b_k)$ we write
$\bar{b}=b_1\cdots b_k$. Let $C_j^k$,  $1\le j\le k$, be the family
of all finite subsets $\sigma=\{\sigma(1),\dots, \sigma(j)\}\subset
\{1,\dots, k\}$ of $j$ different elements. In this case, we write
$\vec{b}_{\sigma}=(b_{\sigma(1)}, \dots, b_{\sigma(j)})$ and
$\bar{b}_{\sigma}=b_{\sigma(1)}\cdots b_{\sigma(j)}$. We also set
$C_0^k=\emptyset$ in which case we understand that
$T_{\vec{b}_{\sigma}}=T$ and $\bar{b}_{\sigma}=1$. If $\sigma\in
C_j^k$ we set $\sigma'=\{1,\dots,k\}\setminus \sigma$ (note that for
$j=0$ we have $\sigma'=\{1,\dots,k\}$). We need the following
multilinear version of \eqref{BMO-M} (see \cite{PT}): given $k\ge 1$,
for any $x\in Q$ we have
\begin{align}\label{m-Holder}
\aver{Q} |f_1\cdots f_k \,h|^{p_0}
&\le
\|f_1\|_{{\rm exp }L}^{p_0}\cdots \|f_k\|_{{\rm exp }L,Q}^{p_0}
\big\| |h|^{p_0}\big\|_{L\,(\log L)^{k\,p_0}}
\nonumber
\\
&\le
\|f_1\|_{{\rm exp }L,Q}^{p_0}\cdots \|f_k\|_{{\rm exp }L,Q}^{p_0}\,
M^{[k\,p_0]+2}\big(|h|^{p_0})(x).
\end{align}

With this in hand and as done with the regular commutators in
Lemma \ref{lemma:comm-apriori} the matter can be reduced to the
case $b_1,\dots,b_k\in L^\infty$. Once we have that, we combine
the ideas from \cite[p. 684]{PT} with the proof above. We write
$F=|T_{\vec{b}} f(x)|^{p_0}\in L^1$ and observe that $F\le
G_Q+H_Q$ where
$$
G_Q= 2^{p_0-1}\,C\,\Big( \sum_{m=1}^{k}\sum_{\sigma\in C_{m}^k}
\overline{(b-\lambda)}_{\sigma}\, |T_{\vec{b}_{\sigma'}} f(x)| +
|T(I-\A_{r(Q)})f_{Q,\vec{b}}(x)|
\Big)^{p_0},
$$
$H_Q=2^{p_0-1}\,|T \A_{r(Q)}f_{Q,\vec{b}}(x)|^{p_0}$, and
$f_{Q,\vec{b}}= \prod_{j=1}^k (b_j-(b_j)_{2\,Q})\,f$. Next, one
estimates $G_Q$, $H_Q$ using the same ideas (with \eqref{m-Holder}
in place of \eqref{BMO-M}):
$$
\aver{Q} G_Q
\lesssim
C\, \Big(\sum_{m=1}^{k}\sum_{\sigma\in C_{m}^k}
M^{[k\,p_0]+2}\big(|T_{\vec{b}_{\sigma'}}
f|^{p_0}\big)(x)+M^{[k\,p_0]+2}\big(|f|^{p_0})(x)\Big)
=
G(x),
$$
and
$$
\Big(\aver{Q} H_Q^q\Big)^{\frac1{q}}
\lesssim
MF(x)+ \sum_{m=1}^{k}\sum_{\sigma\in C_{m}^k}
M^{[k\,p_0]+2}\big(|T_{\vec{b}_{\sigma'}}f|^{p_0}\big)(\overline{x})
\equiv C\,\big(MF(x)+H_2(\overline{x})\big).
$$
>From here the proof proceed as in the case above, noticing that the
length of $\vec{b}_{\sigma'}$ is $k-m\le k-1$ and so the induction
hypothesis applies.
\end{remark}

\subsection{Proof of Theorem \ref{theor:extrapol}}
\label{subsec:proof:theor:extrapol} Assume that the case $p_0=1$ is
proved. Then we show that the general case follows automatically. Set
$\widetilde{p}=p/p_0$, $\widetilde{q}_0=q_0/p_0$ and consider the new
family $\widetilde{\F}$ consisting of the pairs
$(\widetilde{f},\widetilde{g})=(f^{p_0},g^{q_0})$. Observe that $1\le
\widetilde{p}\le \widetilde{q}_0$ and that $\widetilde{p}<\infty$ if
$\widetilde{q}_0=\infty$ (that is, $q_0=\infty$). Besides,
\eqref{extrapol-p} gives that for all
$(\widetilde{f},\widetilde{g})\in \widetilde{\F}$
$$
\int_{\re^n} \widetilde{f}\,^{\widetilde{p}} \,w
\le
C\,\int_{\re^n} \widetilde{g}\, ^{\widetilde{p}} \,w,
\quad
\mbox{for all }w\in A_{\widetilde{p}}\cap
RH_{\left({\widetilde{q}_0}/{\widetilde{p}}\right)'}
$$
provided the left hand side is finite. Therefore, the same holds for
all $1<\widetilde{q}<\widetilde{q}_0$ and  \eqref{extrapol-q} follows
with $q=\widetilde{q}\,p_0$.

Assume now that $p_0=1$. Observe that the case $q_0=\infty$ is
nothing but Rubio de Francia's extrapolation theorem. So we also
impose $q_0<\infty$. The proof of \eqref{extrapol-q} is done on
distinguishing the  two cases $q<p$ and $q>p$. We use the following
notation
$$
\phi(q)=\left(\frac{q_0}{q}\right)'\,(q-1)+1.
$$
Note that $(vii)$ in Proposition \ref{prop:weights} says that if
$q_0/q>1$ then
$
w\in A_{q}\cap RH_{\left(\frac{q_0}{q}\right)'}$ if and only if
$w^{\left(\frac{q_0}{q}\right)'}\in A_{\phi(q)}$. We need the
following auxiliary result  based on Rubio de Francia's algorithm.

\begin{lemma}\label{lemma:algorithms}
Let $1<q<q_0$ and $w$ such that $w\in A_{q}\cap
RH_{\left(\frac{q_0}{q}\right)'}$.
\begin{list}{$(\theenumi)$}{\usecounter{enumi}\leftmargin=.8cm
\labelwidth=0.7cm\itemsep=0.3cm\topsep=0.1cm
\renewcommand{\theenumi}{\alph{enumi}}}

\item If $1\le p<q$ and $0\le h\in L^{(q/p)'}(w)$, then there exists $H\in
L^{(q/p)'}(w)$ such that
\begin{list}{$(\theenumi.\theenumii)\ $}{\usecounter{enumii}\leftmargin=.8cm
\labelwidth=0.7cm\itemsep=0.3cm\topsep=.3cm
\renewcommand{\theenumii}{\arabic{enumii}}}
\item $0\le h\le H$.

\item $\displaystyle \|H\|_{L^{(q/p)'}(w)}\le
2^{\phi(q)'/(q/p)'}\,\|h\|_{L^{(q/p)'}(w)}$.

\item $H\, w\in A_{p}\cap RH_{\left(\frac{q_0}{p}\right)'} $ with constants
independent of $h$.
\end{list}

\item If $q<p\le q_0$ and $0\le h\in L^{q}(w)$, then there exists $H\in
L^q(w)$ such that
\begin{list}{$(\theenumi.\theenumii)\ $}{\usecounter{enumii}\leftmargin=.8cm
\labelwidth=0.7cm\itemsep=0.3cm\topsep=.3cm
\renewcommand{\theenumii}{\arabic{enumii}}}
\item $0\le h\le H$.

\item $\displaystyle \|H\|_{L^{q}(w)}\le
2^{\phi(q)/q}\,\|h\|_{L^{q}(w)}$.

\item $H^{-p/(p/q)'}\, w \in A_{p}\cap RH_{\left(\frac{q_0}{p}\right)'}$ with
constants independent of $h$.
\end{list}
\end{list}
\end{lemma}
Admit  this result for the moment and continue the proof.

\

\noindent \textit{Case $1\le p<q$:} Let $(f,g)\in \F$ be such that
$f,g\in L^q(w)$. Fix $w$ such that $w^{(q_0/q)'}\in A_{\phi(q)}$.
Then,
$$
\|f\|_{L^q(w)}^p
=
\|f^p\|_{L^{q/p}(w)}
=
\sup \int_{\re^n} f^p\,h\,w
$$
where the supremum is taken over all $0\le h \in L^{(q/p)'}(w)$ with
$\|h\|_{L^{(q/p)'}(w)}=1$. Take such a function $h$ and let $H$ be
the corresponding function given by $(a)$ in Lemma
\ref{lemma:algorithms}. Then by $(a.1)$, \eqref{extrapol-p} and
$(a.3)$, we have
$$
\int_{\re^n} f^p\,h\,w
\le
\int_{\re^n} f^p\,H\,w
\le
C\,\int_{\re^n} g^p\,H\,w
$$
provided the middle term is finite. This is indeed the case as by
H\"older's inequality with $q/p>1$ and by $(a.2)$
$$
\int_{\re^n} f^p\,H\,w
\le
\|f\|_{L^q(w)}^{p}\,\|H\|_{L^{(q/p)'}(w)}
\le
2^{\phi(q)'/(q/p)'}\,\|f\|_{L^q(w)}^{p}<\infty.
$$
Note that the same can be done with $g$ and so
$$
\int_{\re^n} g^p\,H\,w
\le
2^{\phi(q)'/(q/p)'}\,\|g\|_{L^q(w)}^{p}.
$$
This readily leads to the desired estimate.

\

\noindent \textit{Case $q<p\le q_0$:} Let $(f,g)\in \F$ be
non-trivial functions such that $f,g\in L^q(w)$. Fix $w$ such that
$w\in A_{q}\cap RH_{\left(\frac{q_0}{q}\right)'}$. We define
$$
h
=
\frac{f}{\|f\|_{L^q(w)}}+\frac{g}{\|g\|_{L^q(w)}}.
$$
Note that $h\in L^q(w)$ and $\|h\|_{L^q(w)}\le 2$. Let $H$ be the
non-negative function given by Lemma \ref{lemma:algorithms} part
$(b)$. Then, using H\"older's inequality with $p/q>1$ we have
\begin{align}
\|f\|_{L^q(w)}
&=
\Big(\,
\int_{\re^n} f^q\,H^{-q/(p/q)'}\,H^{q/(p/q)'}\,w
\Big)^{1/q}
\nonumber
\\
&\le
\Big(\,
\int_{\re^n} f^p\,H^{-p/(p/q)'}\,w
\Big)^{1/p}\,
\Big(\,
\int_{\re^n} H^{q}\,w
\Big)^{\frac1{q\,(p/q)'}}
\nonumber
\\
&\le
C\,\Big(\, \int_{\re^n} f^p\,H^{-p/(p/q)'}\,w
\Big)^{1/p},
\label{case2}
\end{align}
since $(b.2)$ implies
$$
\|H\|_{L^q(w)}
\le
2^{\phi(q)/q}\, \|h\|_{L^q(w)}
\le
2^{1+\phi(q)/q}.
$$
Next, by $(b.1)$ we have $f/\|f\|_{L^q(w)}\le h\le H$. Hence, using
$(b.2)$ we conclude that
\begin{align*}
\Big(\,
\int_{\re^n} f^p\,H^{-p/(p/q)'}\,w
\Big)^{1/p}
& \le
\|f\|_{L^q(w)}\,
\Big(\,
\int_{\re^n} H^{p-p/(p/q)'}\,w
\Big)^{1/p}
=
\|f\|_{L^q(w)}\, \|H\|_{L^q(w)}^{q/p}
\\
&\le
2^{(q/p)(1+\phi(q)/q)}\,\|f\|_{L^q(w)}<\infty.
\end{align*}
This and $(b.3)$ allow us to employ \eqref{extrapol-p}. Hence,
\eqref{case2} yields
$$
\|f\|_{L^q(w)}
\le
C\,
\Big(\,
\int_{\re^n} g^p\,H^{-p/(p/q)'}\,w
\Big)^{1/p}
\le
C\,\|g\|_{L^q(w)}\, \|H\|_{L^q(w)}^{q/p}
\le
C\,\|g\|_{L^q(w)},
$$
where we have used that  $g$ satisfies $g/\|g\|_{L^q(w)}\le H$ due to
$(b.1)$.

To complete the proof it remains to show \eqref{extrapol:v-v}. As in
\cite{CMP} this follows almost automatically from \eqref{extrapol-q}
by changing the family $\F$. Indeed, fix $p_0<r<q_0$ and given
$\{(f_j,g_j)\}_j\subset\F$ we define
$$
F_r
=
\Big(\sum_j f_j^r\Big)^{1/r},
\qquad\qquad
G_r
=
\Big(\sum_j g_j^r\Big)^{1/r}.
$$
We consider a new family $\F_r$ consisting of all the pairs
$(F_r,G_r)$. Observe that if $(F_r,G_r)\in \F_r$, using
\eqref{extrapol-q} with $q=r$,  we have

$$
\|F_r\|_{L^r(w)}^r
=
\sum_j\int_{\re^n} f_j^r\,w
\le
C\,\sum_j\int_{\re^n} g_j^r\,w
=
C\,\|G_r\|_{L^r(w)}^r,
$$
for all $w\in A_{r/p_0}\cap RH_{(q_0/r)'}$. This means that the
family $\F_r$ satisfies \eqref{extrapol-p} with $p=r$. Thus, as we
have just obtained, it satisfies \eqref{extrapol-q} for all
$p_0<q<q_0$ which turns out to be \eqref{extrapol:v-v}. \qed

\begin{proof}[Proof of Lemma \ref{lemma:algorithms}]
We first observe that
$$
w^{(q_0/q)'}\in A_{\phi(q)}
\qquad
\Longleftrightarrow
\qquad
w^{1-q'}=w^{(q_0/q)'(1-\phi(q)')}\in A_{\phi(q)'}.
$$
Given any weight $0<u<\infty$ a.e. we define the operator
$$
S_u f
=
\frac{M(f\,u)}{u}.
$$
This operator will be used to perform different versions of Rubio de
Francia's algorithm. We start with $(a)$: let $1\le p<q$ and $h\in
L^{(q/p)'}(w)$. We set $u=w^{q'/\phi(q)'}$. Then, as $w^{1-q'}\in
A_{\phi(q)'}$ we have
\begin{align*}
\|S_u f\|_{L^{\phi(q)'}(w)}^{\phi(q)'}
&=
\int_{\re^n} M(f\,u)^{\phi(q)'}\,u^{-\phi(q)'}\,w
= \int_{\re^n}
M(f\,u)^{\phi(q)'}\, w^{1-q'}
\\
&\le
C\,\int_{\re^n} |f\,u|^{\phi(q)'}\, w^{1-q'}
=
C\, \|f\|_{L^{\phi(q)'}(w)}^{\phi(q)'}.
\end{align*}
Let us write $\|S_u\|$ for the norm of $S_u$ as a bounded operator on
$L^{\phi(q)'}(w)$. We define the following version of Rubio de
Francia's algorithm: for $0\le f\in L^{\phi(q)'}(w)$
$$
\R f
=
\sum_{k=0}^{\infty} \frac{S_u^k f}{2^k\,\|S_u\|^k},
$$
where $S_u^k$ is the $k$-iteration of the operator $S_u$ for $k\ge 1$
and  $S_u^0$ is the identity operator. Given $0\le h\in
L^{(q/p)'}(w)$ we define
$$
H
=
\R\big(\,h^{(q/p)'/\phi(q)'}\big)^{\phi(q)'/(q/p)'}.
$$
Note that
$$
0\le f\le \R f,
\qquad\qquad
\|\R f\|_{L^{\phi(q)'}(w)}\le 2\,\|f\|_{L^{\phi(q)'}(w)},
$$
and so $H$ satisfies $(a.1)$ and $(a.2)$. Note that we also have
$$
S_u(\R f)\le 2\,\|S_u\|\, \R f
\
\Longleftrightarrow
\
M(u\, \R f )
\le
C\, u \,\R f
\
\Longleftrightarrow
\
u\, \R f\in A_1
$$
and therefore $H^{(q/p)'/\phi(q)'}\,u\in A_1$ with constant
independent of $h$. Then for all cube $Q\subset \re^n$ (the averages
are with respect to Lebesgue measure)
\begin{equation}\label{A1:case1}
\aver{Q} H^{(q/p)'/\phi(q)'}\,u
\le
C\,H^{(q/p)'/\phi(q)'}(x)\,u(x),
\qquad\mbox{a.e. }x\in Q.
\end{equation}
We show $(a.3)$, that is, $(H\,w)^{(q_0/p)'}\in A_{\phi(p)}$. If
$p=1$ then \eqref{A1:case1} turns out to be
$$
\aver{Q} \big(H\,w\big)^{q_0'}
\le
C\,\big(H(x)\,w(x)\big)^{q_0'},
\qquad\mbox{a.e. }x\in Q,
$$
that is, $(H\,w)^{q_0'}\in A_{\phi(1)}=A_1$ as desired. If $p>1$,
using \eqref{A1:case1} we have
\begin{align*}
I
&=
\aver{Q} (H\,w)^{(q_0/p)'\,(1-\phi(p)')}
=
\aver{Q} (H\,w)^{1-p'}
\\
&\lesssim
\Big(
\aver{Q} H^{(q/p)'/\phi(q)'}\,u
\Big)^{-\frac{(p'-1)\,\phi(q)'}{(q/p)'}}
\,
\Big(\aver{Q} u^{\frac{(p'-1)\phi(q)'}{(q/p)'}}\,
w^{1-p'}\Big)
\\
&=
\Big(
\aver{Q} H^{\frac{q_0\,(q-1)}{(q_0-1)\,(q-p)}}\,u
\Big)^{-\frac{(q_0-1)\,(q-p)}{q_0\,(q-1)\,(p-1)}}
\Big(\aver{Q}
w^{1-q'}\Big)
=
I_1\cdot I_2.
\end{align*}
Since $1<p<q<q_0$ we have that
$$
s
=
\frac{q_0\,(q-1)}{(q_0-1)\,(q-p)}\,\frac1{(q_0/p)'}
=
\frac{(q-1)\,(q_0-p)}{(q_0-1)\,(q-p)}>1,
\qquad\quad
s'
=
\frac{(q-1)\,(q_0-p)}{(p-1)\,(q_0-q)} \ .
$$
Then by H\"older's inequality we obtain
\begin{align*}
II
&=
\aver{Q} \big(\,H\,w\big)^{(q_0/p)'}
\le
\Big(\aver{Q} H^{(q_0/p)'\,s}\,u\Big)^{1/s} \,
\Big(\aver{Q}
w^{(q_0/p)'\,s'}\,u^{1-s'}\Big)^{1/s'}
\\
&=
\Big(\aver{Q}
H^{\frac{q_0\,(q-1)}{(q_0-1)\,(q-p)}}\,u\Big)^{1/s} \,
\Big(\aver{Q}
w^{(q_0/q)'}\Big)^{1/s'}
=
II_1\cdot II_2.
\end{align*}
We gather $I_1$ and $II_1$:
$$
I_1^{\phi(p)-1}\cdot II_1
=
\Big(
\aver{Q} H^{\frac{q_0\,(q-1)}{(q_0-1)\,(q-p)}}\,u
\Big)^{\frac{1}{s}-(\phi(p)-1)\frac{(q_0-1)\,(q-p)}{q_0\,(q-1)\,(p-1)}}
=1
$$
since the outer exponent is equal to $0$. On the other hand, for
$I_2$ and $II_2$ we observe that
\begin{eqnarray*}
I_2^{\phi(p)-1}\cdot II_2
&=&
\Big(\aver{Q} w^{1-q'}\Big)^{\phi(p)-1}
\,
\Big(\aver{Q}
w^{(q_0/q)'}\Big)^{1/s'}
\\
&=&
\bigg[
\Big(\aver{Q}
w^{(q_0/q)'\,(1-\phi(q)')}\Big)^{(\phi(p)-1)\,s'}
\Big(\aver{Q}
w^{(q_0/q)'}\Big)
\bigg]^{1/s'}
\\
&=&
\bigg[
\Big(\aver{Q}
w^{(q_0/q)'}\Big) \,
\Big(\aver{Q}
w^{(q_0/q)'\,(1-\phi(q)')}\Big)^{\phi(q)-1}
\bigg]^{1/s'}
\le
C,
\end{eqnarray*}
since $w^{(q_0/q)'}\in A_{\phi(q)}$. As a consequence of these latter
estimates we can conclude that $(H\,w)^{(q_0/p)'}\in A_{\phi(p)}$:
\begin{align*}
\lefteqn{\hskip-2cm
\Big(\aver{Q} \big(\,H\,w\big)^{(q_0/p)'}\Big)\,
\Big(\aver{Q}
(H\,w)^{(q_0/p)'\,(1-\phi(p)')}\Big)^{\phi(p)-1}
=
I^{\phi(p)-1}\cdot II }
\\
&\le
C\,(I_1^{\phi(p)-1}\cdot II_1)\, (I_2^{\phi(p)-1}\cdot II_2)\,
\le C .
\end{align*}

We now prove $(b)$.  Let $h\in L^q(w)$ and
$u=w^{(1-(q_0/q)')/\phi(q)}$. Since $w^{(q_0/q)'}\in A_{\phi(q)}$ we
have
\begin{align*}
\|S_u f\|_{L^{\phi(q)}(w)}^{\phi(q)}
&=
\int_{\re^n} M(f\,u)^{\phi(q)}\,u^{-\phi(q)}\,w = \int_{\re^n}
M(f\,u)^{\phi(q)}\, w^{(q_0/q)'}
\\
& \le
C\,\int_{\re^n} |f\,u|^{\phi(q)}\, w^{(q_0/q)'}
=
C\, \|f\|_{L^{\phi(q)}(w)}^{\phi(q)}.
\end{align*}
Let us write $\|S_u\|$ for the norm of $S_u$ as a bounded operator on
$L^{\phi(q)}(w)$. Rubio de Francia's algorithm to be used now is
given by
$$
\R f
=
\sum_{k=0}^{\infty} \frac{S_u^k f}{2^k\,\|S_u\|^k},
$$
for $0\le f\in L^{\phi(q)}(w)$. Given $0\le h\in L^{q}(w)$ we define
$$
H
=
\R\big(\,h^{q/\phi(q)}\big)^{\phi(q)/q}.
$$
Note that
$$
0\le f\le \R f,
\qquad\qquad
\|\R f\|_{L^{\phi(q)}(w)}\le 2\,\|f\|_{L^{\phi(q)}(w)},
$$
and so $H$ satisfies $(b.1)$ and $(b.2)$. As in the other case
$$
S_u(\R f)\le 2\,\|S_u\|\, \R f
\
\Longleftrightarrow
\
M(u\, \R f )
\le
C\, u\, \R f
\
\Longleftrightarrow
\
u\, \R f\in A_1
$$
and so $H^{q/\phi(q)}\,u\in A_1$ with constant independent of $h$.
Thus for all cubes $Q\subset \re^n$
\begin{equation}\label{A1:case2}
\aver{Q} H^{q/\phi(q)}\,u
\le
C\, H^{q/\phi(q)}(x)\,u(x),
\qquad\mbox{a.e. }x\in Q.
\end{equation}
We prove $(b.3)$. We do first the case $p=q_0$ and we have to see
that $H^{-(q_0-q)}\,w\in A_{q_0}\cap RH_\infty$. Note that
\eqref{A1:case2} can be rewritten as
$$
\aver{Q} \big(H^{q_0-q}\,w^{-1})^{q_0'-1}
\le
C\, \big(H^{q_0-q}(x)\,w^{-1}(x))^{q_0'-1},
\qquad\mbox{a.e. }x\in Q.
$$
Then, for almost every $x\in Q$ we have
$$
H^{-(q_0-q)}(x)\,w(x)
\lesssim
\Big(\aver{Q}
\big(H^{q_0-q}\,w^{-1})^{q_0'-1}\Big)^{-\frac1{q_0'-1}}
\le
\aver{Q}  H^{-(q_0-q)}\,w
$$
where in the last estimate we have used Jensen's inequality with the
convex function $t\mapsto t^{-1/(q_0'-1)}$. This shows that
$H^{-(q_0-q)}\,w\in RH_\infty$. On the other hand, we also have
$$
\Big(\aver{Q}  H^{-(q_0-q)}\,w\Big)
\lesssim
\Big(\aver{Q}
\big(H^{-(q_0-q)}\,w\big)^{1-q_0'}\Big)^{-(q_0-1)}
$$
which automatically implies that $H^{-(q_0-q)}\,w\in A_{q_0}$. This
completes the case $p=q_0$.

If $p<q_0$, $(b.3)$ is equivalent to $\big(\,H^{-p/(p/q)'}\,
w\big)^{(q_0/p)'}\in A_{\phi(p)}$. By \eqref{A1:case2} we observe
that
\begin{align*}
I
&=
\aver{Q} \big(\,H^{-p/(p/q)'}\, w\big)^{(q_0/p)'}
\lesssim
\Big(
\aver{Q}  H^{q/\phi(q)}\, u
\Big)^{-\frac{p\,(q_0/p)'\,\phi(q)}{(p/q)'\,q}}
\,
\Big(
\aver{Q}  u^{\frac{p\,(q_0/p)'\,\phi(q)}{(p/q)'\,q}}\, w^{(q_0/p)'}
\Big)
\\
&=
\Big(
\aver{Q}  H^{q/\phi(q)}\, u
\Big)^{-\frac{p\,(q_0/p)'\,\phi(q)}{(p/q)'\,q}}
\,
\Big(
\aver{Q}  w^{(q_0/q)'}
\Big)
=
I_1\cdot I_2.
\end{align*}
Since $1<q<p<q_0$ we have that
$$
s
=
\frac{q\,(p-1)}{\phi(q)\,(p-q)}
=
\frac{(q_0-q)\,(p-1)}{(q_0-1)\,(p-q)}
>1,
\qquad\quad
s'
=
\frac{(q_0-q)\,(p-1)}{(q_0-p)\,(q-1)}.
$$
By H\"older's inequality we obtain
\begin{align*}
II
&=
\aver{Q}  \big(\,H^{-p/(p/q)'}\, w\big)^{(q_0/p)'\,(1-\phi(p)')}
=
\aver{Q}  \big(\,H^{-p/(p/q)'}\, w\big)^{1-p'}
=
\aver{Q}  H^{\frac{p-q}{p-1}}\, w^{1-p'}
\\
&\le
\Big(\aver{Q}  H^{\frac{p-q}{p-1}\,s}\,u \Big)^{1/s} \,
\Big(\aver{Q}  w^{(1-p')\,s'}\,u^{1-s'} \Big)^{1/s'}
=
\Big(
\aver{Q}  H^{q/\phi(q)}\,u
\Big)^{1/s} \,
\Big(\aver{Q}  w^{1-q'}\Big)^{1/s'}
\\
&=
II_1\cdot II_2.
\end{align*}
For $I_1$ and $II_1$ we have
$$
I_1\cdot II_1^{\phi(p)-1}
=
\Big(
\aver{Q}  H^{q/\phi(q)}\, u
\Big)^{-\frac{p\,(q_0/p)'\,\phi(q)}{(p/q)'\,q}+\frac{\phi(p)-1}{s}}
=1
$$
since the outer exponent vanishes. On the other hand, since
$w^{(q_0/q)'}\in A_{\phi(q)}$,
\begin{align*}
I_2\cdot II_2^{\phi(p)-1}
&=
\Big(
\aver{Q}  w^{(q_0/q)'}
\Big)
\,\Big(\aver{Q} w^{1-q'}\Big)^{\frac{\phi(p)-1}{s'}}
\\
&=
\Big(
\aver{Q}  w^{(q_0/q)'}
\Big)
\,\Big(\aver{Q} w^{(q_0/q)'(1-\phi(q)')}\Big)^{\phi(q)-1}
\le C.
\end{align*}
Collecting the last two  estimates we conclude that
$\big(H^{-p/(p/q)'}\, w\big)^{(q_0/p)'}\in A_{\phi(p)}$:
\begin{align*}
\lefteqn{\hskip-2.8cm
\Big(\aver{Q}  \big(\,H^{-\frac{p}{(p/q)'}}
w\big)^{(q_0/p)'}\Big) \,
\Big(
\aver{Q}  \big(\,H^{-\frac{p}{(p/q)'}}
w\big)^{(q_0/p)'\,(1-\phi(p)')}
\Big)^{\phi(p)-1}
=
I\cdot II^{\phi(p)-1} }
\\
&
\le C\,
(I_1\cdot II_1^{\phi(p)-1})\, (I_2\cdot II_2^{\phi(p)-1})
\le
C.
\end{align*}
\end{proof}

\part{Calder\'on-Zygmund methods}\label{part:two}

\section{Introduction}

This section develops a circle of ideas based on the Calder\'on-Zygmund
decomposition. This decomposition was invented in the celebrated
article \cite{CZ} to prove that certain singular integrals of
convolution type are of weak-type $(1,1)$. Recall that this
decomposition is non-linear and breaks up $L^1$ functions into good
and bad parts. The good part is bounded, while the bad part is a sum
of localized and oscillating functions. The oscillation is in the
sense of a vanishing mean.  This turned out to be a very versatile
tool.

The application towards singular integrals was refined in \cite{Hor}
with a minimal regularity condition on the kernel matching the
oscillation of the bad parts. Then, this was generalized to what is
now called  Calder\'on-Zygmund operators, see, e.g, \cite{Mey}. We note
that a key ingredient in these arguments is the \textit{a priori}
strong or weak-type $(p_{0}, p_{0})$ of the operator for some
$p_{0}>1$.

Kernel regularity  in some sense is needed for such arguments. After
the results obtained in \cite{Heb} and \cite{DR} in a functional
calculus setting, a general weak-type  (1,1) criterion is formulated
in \cite{DMc}.  It still exploits the Calder\'on-Zygmund decomposition
but does not use the oscillation of the bad part. The regularity is
expressed in the integrability properties of the kernel of
$T(Id-\A_{r})$ where $\A_{r}$ , $r>0$, is some approximation to the
identity. In the classical case, $\A_{r}$ would be an ordinary
mollifying operator with a smooth bump function.

\cite{BK1} develops this idea further for  singular ``non-integral''  operators and
establishes a weak-type $(p,p)$ criterion,  still assuming of course
\textit{a priori}  weak-type $(p_{0}, p_{0})$ boundedness for some
$p_{0}>p$. This result is presented in \cite{Aus-mem} with a simpler
and stronger statement.  This is typically an unweighted result but
as it works in spaces of homogeneous type, it applies with underlying
doubling measure $w(x) \, dx$,  $w\in A_{\infty}$.

In a sense, we have not much to add to this story. However we present
it once again as its argument is needed for  further development
(Section \ref{section:BK}). First, a slight strengthening of the
hypotheses yields for free boundedness results for  commutators of
the operator with bounded mean oscillation functions (Section
\ref{section:BK:comm}). Second, we observe that similar unweighted
estimates plus an \textit{a priori} weighted weak-type
$(p_{0},p_{0})$ estimate of $T$ implies  weighted weak-type $(p,p)$
estimate  for  a range  of $p$'s with $p<p_{0}$ depending on the
class of weights (Section \ref{section:BK:weights}).

We also present in  Section \ref{section:special-CZ}  a result of
independent interest but needed in \cite{AM3} concerning  a
Calder\'on-Zygmund decomposition for a function in $\RR^n$ with
gradient controlled in some $L^p(w)$ space for some $p\ge 1$ and
doubling weight $w$ supporting a Poincar\'e inequality. Such  a
decomposition is used is \cite{Aus-pubm} in the Euclidean setting and
a similar decomposition appear earlier in
\cite{CM} and \cite{BS} for the purpose of real interpolation for Sobolev spaces.
See also \cite{AC} for an extension to Riemannian manifolds.

\section{Extended Calder\'on-Zygmund theory}\label{section:CZ}

Except for Section \ref{section:BK:SHT}, we work in $\re^n$ endowed with a Borel
doubling measure $\mu$ (and we remind the reader that in applications
$d\mu(x)=w(x)\,dx$ with $w\in A_\infty$).

\subsection{Blunck and Kunstmann's theorem}\label{section:BK}
We use the following notation: if $B$ is a ball with radius $r(B)$
and $\lambda>0$, $\lambda\, B$ denotes the concentric ball with
radius $r(\lambda \, B)= \lambda\, r(B)$, $C_{j}(B)=2^{j+1}\,
B\setminus 2^j\, B$ when $j\ge 2$, $C_{1}(B)=4B$, and
\begin{equation}\label{eq:meanCj}
\aver{C_j(B)} h\,d\mu
=
\frac1{\mu(2^{j+1}B)}\,\int_{C_{j}(B)} h\,d\mu.
\end{equation}
We say that the doubling measure $\mu$ has doubling order $D>0$ if
$\mu(\lambda\,B)\le C_\mu\,\lambda^D\,\mu(B)$ for every ball $B$
and every $\lambda>0$.

The following result appears in a paper by Blunck and Kunstmann
\cite{BK1} in a slightly more complicated way with extra hypotheses. This version is due to  one of us \cite{Aus-mem}.

\begin{theor}\label{theor:B-K:small}
Let $\mu$ be a doubling Borel measure on $\re^n$ with doubling order
$D$ and $1\le p_0<q_0\le \infty$. Suppose that $T$ is a sublinear
operator of weak-type $(q_0,q_0)$. Let $\D$ be a subspace of
$L^{q_0}(\mu)\cap L^{p_{0}}(\mu)$ stable under truncation by indicator
functions of measurable sets. Let $\{\A_r\}_{r>0}$ be a family of
operators acting from $\D$ into $L^{q_0}(\mu)$. Assume that for $j\ge
2$,
\begin{equation}\label{small:T:I-A}
\Big( \aver{C_j(B)} |T(I-\A_{r(B)})f|\,d\mu\Big)
\le
\alpha_j\,\Big(\aver{B} |f|^{p_0}\,d\mu\Big)^{\frac1{p_0}}
\end{equation}
and for $j\ge 1$
\begin{equation}\label{small:A}
\Big(
\aver{C_j(B)} |\A_{r(B)}f|^{q_0}\,d\mu\Big)^{\frac1{q_0}}
\le
\alpha_j\,\Big(\aver{B} |f|^{p_0}\,d\mu\Big)^{\frac1{p_0}},
\end{equation}
for all ball $B$ with $r(B)$ its radius and for all $f\in \D$
supported in $B$. If $\sum_j \alpha_j\,2^{D\,j}<\infty$ then $T$ is
of weak-type $(p_0,p_0)$ and hence $T$ is of strong-type $(p,p)$ for
all $p_0<p<q_0$. More precisely, there exists a constant $C$ such
that for all $f \in \D$,
$$
\|T f\|_{L^p(\mu)} \le C\,  \| f\|_{L^p(\mu)}.
$$
\end{theor}

\medskip

\subsection{Commutators with BMO functions: part II}\label{section:BK:comm}

A slightly strengthening of the hypotheses above yields an analog result
for the commutators with bounded mean oscillation functions. In this
case, since the underlying measure is $\mu$, we work with functions
$b\in \BMO(\mu)$ (the definition is as the classical one replacing
$dx$ by $\mu$). As $\mu$ is a doubling measure, John-Nirenberg's
inequality holds in $\BMO(\mu)$. The
definition of the commutator is the same as in Section \ref{section:comm:p-big} but in this case we assume
that $T$ is of weak-type $(q_0,q_0)$ in place of being bounded on
$L^{p_0}$. This  still guarantees that the commutator is well defined.

\begin{theor}\label{theor:small:comm}
Let $\mu$ be a doubling Borel measure on $\re^n$ with doubling order
$D$,  $1\le p_0<q_0\le \infty$, $b\in \BMO(\mu)$ and  $k\in \NN$,
$k\ge 1$. Suppose that $T$ is a sublinear operator and that $T$ and
$T_b^m$ for $m=1,\dots,k$ are of weak-type $(q_0,q_0)$. Let
$\{\A_r\}_{r>0}$ be a family of operators acting from
$L^\infty_c(\mu)$ into $L^{q_0}(\mu)$. Assume that for any ball $B$
with $r(B)$ its radius and for all $f\in L^\infty_c$ supported in
$B$, \eqref{small:A} holds, and \eqref{small:T:I-A} is replaced by
the stronger assumption
\begin{equation}\label{small:T:I-A:comm}
\Big( \aver{C_j(B)} |T(I-\A_{r(B)})f|^{r}\,d\mu\Big)^{\frac1{r}}
\le
\alpha_j\,\Big(\aver{B} |f|^{p_0}\,d\mu\Big)^{\frac1{p_0}}
\end{equation}
for some $r>1$ and all $j\ge 2$. If $\sum_j
\alpha_j\,2^{D\,j}\,j^k<\infty$ then for all $p_0<p<q_0$, there
exists a constant $C$ \textup{(}independent of $b$\textup{)} such
that for all $f \in L^\infty_c(\mu)$,
$$
\|T_b^k f\|_{L^p(\mu)} \le C\,  \|b\|_{\BMO(\mu)}^k\, \|
f\|_{L^p(\mu)}.
$$
\end{theor}

\begin{remark}\label{remark:indu-0:small}\rm
Under the  assumptions above,  we have  $\sum_j
\alpha_{j}\, \,2^{D\,j}<\infty$ and consequently, Theorem
\ref{theor:B-K:small} implies that $T=T_b^0$ is of weak-type
$(p_0,p_0)$ and hence bounded on $L^p(\mu)$ for all $p_0<p<q_0$.
\end{remark}

\begin{remark}\rm
 In applications we will use this result with underlying
measure $d\mu(x)=w(x)\,dx$ with $w\in A_\infty$ and so the weight is
hidden in the measure. Let us mention that if $w\in
A_\infty$, and so $dw$ is a doubling measure, then the reverse H\"older
property yields that $\BMO(w)=\BMO$ with equivalent norms.
\end{remark}

\begin{remark}\rm
Our argument requires that the commutators are already weak-type
$(q_{0},q_{0})$, which could make this result useless. However, this
hypothesis can be obtained from Theorem \ref{theor:comm}, see
\cite{AM3} for examples of this.
\end{remark}

\begin{remark}\label{remark:multi-comm:p-small}\rm
As in Remark \ref{remark:multi-comm:p-big}, we can also consider
multilinear commutators associated with a vector of symbols
$\vec{b}=(b_1,\dots,b_k)$ with entries in $\BMO(\mu)$. In this case,
we can formulate an analog of Theorem \ref{theor:small:comm} proving
that $T_{\vec{b}}$ is bounded on $L^p(\mu)$ (see Remark
\ref{remark:proof:multi-comm:p-small} below). The precise statement
is left to the reader.
\end{remark}

\subsection{Weighted estimates}\label{section:BK:weights}

We present the following weighted version of Theorem
\ref{theor:B-K:small} which is used in \cite{AM4}.

\begin{theor}\label{theor:B-K:small:w}
Let $\mu$ be a doubling Borel measure on $\re^n$, $w\in A_{\infty}$
with doubling order $D_{w}$. Let $\D_{1} \subset \D_{2}$ be subspaces
of $L^{q_{0}}(w)$ and suppose that they are
 stable under truncation by indicator
functions of measurable sets. Let $T$ be a sublinear operator defined
on    $\D_{2}$. Let $\{\A_r\}_{r>0}$ be a family of operators acting
from $\D_{1}$ into $\D_{2}$. Let  $1\le p_0<q_0\le \infty$. Assume the following conditions
\begin{list}{$(\theenumi)$}{\usecounter{enumi}\leftmargin=1cm
\labelwidth=0.7cm\itemsep=0.1cm\topsep=.3cm
\renewcommand{\theenumi}{\alph{enumi}}}

\item There exists  $q\in \W_{w}(p_{0},q_{0})$ such that  $T$ is bounded from $L^q(w)$ to $L^{q,\infty}(w)$.

\item For all $j\ge 1$, there exist constants $\alpha_{j}$ such that   for any ball $B$ with $r(B)$ its radius and
for any $f\in \D_{1}$ supported in $B$,
\begin{equation}\label{small:A:w}
\Big(
\aver{C_j(B)} |\A_{r(B)}f|^{q_0}\,d\mu\Big)^{\frac1{q_0}}
\le
\alpha_j\,\Big(\aver{B} |f|^{p_0}\,d\mu\Big)^{\frac1{p_0}}.
\end{equation}

\item There exists $\beta>(s_{w})'$, \textit{i.e.} $w\in RH_{\beta'}$,
with the following property: for all $j\ge 2$, there exist constants
$\alpha_{j}$ such that  for any ball $B$ with $r(B)$ its radius and
for any $f\in \D_{1}$ supported in $B$ and for  $j\ge 2$,
\begin{equation}\label{small:T:I-A:w}
\Big( \aver{C_j(B)} |T(I-\A_{r(B)})f|^{\beta}\,d\mu\Big)^{1/\beta}
\le
\alpha_j\,\Big(\aver{B} |f|^{p_0}\,d\mu\Big)^{\frac1{p_0}}.
\end{equation}
\item $\sum_j \alpha_j\,2^{D_{w}\,j}<\infty$ for $\alpha_{j}$ in (b) and (c).
\end{list}
   Then $T$ is of strong-type
$(p,p)$ with respect to  $w$  for all $p\in \W_{w}(p_{0},q_{0})$ with
$p<q$. More precisely,  for such a $p$, there exists a constant $C$ such that for all
$f \in \D_{1}$,
$$
\|T f\|_{L^p(w)} \le C\,  \| f\|_{L^p(w)}.
$$
\end{theor}

\begin{proof}
Fix a ball $B$, $f$ supported in $B$ and let $g=|T(I-\A_{r(B)})f|$
and $h=|\A_{r(B)}f|$. Let $p\in \W_{w}(p_{0},q_{0})$ with $p<q$.
Since $w\in RH_{(q_{0}/q)'}$ and $w\in A_{p/p_{0}}$,
\eqref{small:A:w} yields
\begin{equation*}
\Big(
\aver{C_j(B)} h^{q}\,dw\Big)^{\frac1{q}}
\lesssim
\Big(
\aver{C_j(B)} h^{q_0}\,d\mu\Big)^{\frac1{q_0}}
\le
\alpha_j\,\Big(\aver{B} |f|^{p_0}\,d\mu\Big)^{\frac1{p_0}}
\lesssim
\alpha_j\,\Big(\aver{B} |f|^{p}\,dw\Big)^{\frac1{p}}.
\end{equation*}
Then as
$w\in RH_{\beta'}$ and $w\in A_{p/p_{0}}$, \eqref{small:T:I-A:w}
implies
\begin{equation*}
\aver{C_j(B)} g \,dw
\lesssim \Big( \aver{C_j(B)} g^{\beta}\,d\mu\Big)^{1/\beta}
\le
\alpha_j\,\Big(\aver{B} |f|^{p_0}\,d\mu\Big)^{\frac1{p_0}}
\lesssim \alpha_j\,\Big(\aver{B} |f|^{p}\,dw\Big)^{\frac1{p}}.
\end{equation*}
Thus we are back to the hypothesis of Theorem
\ref{theor:B-K:small} for the
doubling measure $w\,d\mu$ and with exponents $p<q$.  This implies
that $T$ has weak-type $(p,p)$ with respect to  $w\, d\mu$. As $p$
is arbitrary in an open interval, this implies also strong-type by
Marcinkiewicz interpolation theorem.
\end{proof}

\begin{remark}\rm
Note that  \eqref{small:A:w} and \eqref{small:T:I-A:w}  are unweighted
assumptions. Since we assume weighted weak-type $(q,q)$ for $T$, this
seems useless in applications. In fact, it is a good companion of
Theorem \ref{theor:main-w}. See the application to Riesz transforms
on manifolds in \cite{AM4}.
\end{remark}

\begin{remark}\rm
An examination of the argument shows that if in addition $w\in A_{1}$ then
weighted weak-type holds at $p=p_0$.
\end{remark}

\begin{remark}\rm
A simple and special case is the following. If $(b)$, $(c)$ and $(d)$ hold for $p_{0}=1$ and $q_{0}=\infty$,  then it suffices that $(a)$ holds for some $q$ with $q>r_{w}$ and the conclusion
holds for all $p\in (r_{w},q)$.
\end{remark}

\begin{remark}\rm
We can obtain a version of Theorem \ref{theor:B-K:small:w}  for
commutators with BMO functions: let $k\ge 1$, $b\in \BMO$ and $w\in
A_\infty$. In $(a)$ we further assume that $T_b^m$, for
$m=1,\dots,k$, are bounded from $L^q(w)$ to $L^{q,\infty}(w)$; the
series in $(d)$ becomes $\sum_j
\alpha_j\,2^{D_w\,j}\,j^k<\infty$; $(b)$, $(c)$ remain the same. In
such a case, we show that $T_b^k$ is bounded on $L^p(w)$ for $p<q$,
$p\in \W_w(p_0,q_0)$.

The proof is almost identical and we only give the main ideas. The
computations for $h$ do not change. To estimate $g$, in the left-hand
side, we need to start with an $L^r(w)$-norm in place of the
$L^1(w)$-norm. We pick $r>1$ so that $(s_w)'<\beta/r<\beta$ (note
that $(s_w)'<\beta$). This guarantees that $w\in RH_{(\beta/r)'}$ and
from the $L^r(w)$-norm we pass to  the $L^\beta(\mu)$-norm, after
this the desired estimate follows in the same manner. Thus, we can
apply Theorem \ref{theor:small:comm} to obtain that $T_b^k$ is
bounded on $L^{\tilde{p}}(w)$ for all $p<\tilde{p}<q$. As $p$ is
arbitrary in an open interval, we conclude that $T_b^k$ is bounded on
$L^{p}(w)$ for all $p<q$ such that $p\in \W_w(p_0,q_0)$.
\end{remark}

\subsection{Extension to spaces of homogeneous type}\label{section:BK:SHT}

The preceding results in this part have been obtained in $\re^n$ equipped with
a doubling measure $\mu$. In \cite{AM3} we will use them with
$\mu$ being either the Lebesgue measure or $d\mu(x)=w(x)\,dx$ with
$w\in A_\infty$ and in \cite{AM4}, $\RR^n$ will be replaced by a manifold or a Lie group.  It is not difficult to see that all the proofs
can be adapted to the case of general spaces of homogeneous type
$(\X,d,\mu)$ (see \cite{CW}, \cite{Chr}, \cite{Ste}). Precise
statements and details are left to the reader.

Let us just make a point about the definition \eqref{eq:meanCj}. It
would have looked more natural to  use the ``true'' mean of $h$
over $C_j(B)$ where we divide by $\mu(C_j(B))$ in place of
$\mu(2^{j+1}\,B)$. Our choice is justified partly by the fact that we do not know whether
 $2\, B\setminus B$ and $2\, B$ have comparable mass for all balls, and partly  since (fortunately) $\mu(2^{j+1}\,B)$ is the quantity that appears in computations. Let us note a fairly weak sufficient condition on $\X$ insuring this comparability (which is surely known but we could not find an explicit statement in the literature)

\begin{lemma} Assume that there exists  $\ep\in (0,1)$ such that for any ball $B\subset \X$,  $(2-\ep)\, B \setminus B \ne \emptyset$. Then, $\mu(2\,B\setminus B)\approx \mu(2\,B)$ for any ball $B$, where the implicit constants are independent of $B$.
\end{lemma}

It would be nice to be able to take $\ep=0$ in the above statement. The argument below  shows that $\mu(2\,B\setminus B) \ge C \mu(2\,B)$ but with $C$ depending on $B$. So our statement is the next best thing.

We prove the lemma. It suffices to show that $\mu(2\,B) \ge \nu\,
\mu(B)$ for some $\nu >1$. Choose $1<c< \frac{3}{3-\ep}$. Let $B$ be
a ball, $x_B$ its center and $r$ its radius. By hypothesis, there
exists $x\in B(x_B,(2-\ep)\,c\,r)\setminus B(x_B,c\,r)$. Set
$B'=B(x,(c-1)\,r)$ and note that $B'\subset 2\,B\setminus B$. Thus
$\mu(2B)\ge \mu(B) + \mu(B')$. Now $B \subset \kappa\,B'$ with
$\kappa= \frac{(3-\ep)c}{c-1}$, hence $\mu(B)
\le \mu(\kappa\,B') \le C \kappa^D \mu(B')$ where $D$ is the doubling order of $\mu$. Therefore, $\mu(2\,B)\ge (1+ (C \kappa^D)^{-1})\,\mu(B)$ as desired.

\

Remark that if we had  assumed that  \textit{all} annuli are non-empty then we would obtain  for \textit{all} $\lambda>1$,
$\mu(\lambda\,B)\ge c_\mu\,\lambda^{d}\,\mu(B)$ for some $c_\mu\ge 1$
and $d>0$ depending on $\mu$.  Let us  finally
observe that Theorems \ref{theor:B-K:small}, \ref{theor:small:comm}
and \ref{theor:B-K:small:w} hold  with $a$-adic annuli for
some fixed $a>1$ instead of dyadic ones. The needed changes in the statements and proofs are left
to the reader.

\section{On a special Calder\'on-Zygmund decomposition}\label{section:special-CZ}

The standard Calder\'on-Zygmund decomposition of functions allows one
to decompose a function into a sum of a good bounded function and bad
but localized functions.  This decomposition depends on the level
sets of the maximal function of $f$. This is used to prove
boundedness results such as  Theorem \ref{theor:B-K:small}.

If one wants to prove estimates like $\|Tf\|_{p} \lesssim
\sum_{j=1}^n \|\partial_{j} f\|_{p}$  then one observes that   the
level sets under control are those of the maximal function of each
partial $\partial_{j}f$. But unless one can explicitly express $Tf$
in terms of the functions $\partial_{j}f$, the decomposition applied
to each $\partial_{j}f$ does not allow to split $f$ as before.

The idea of the following lemma, which is applied in \cite{AM3}, is
to split $f$ according to some information on its gradient. This was
done in \cite{Aus-mem} for Lebesgue measure in $\RR^n$. We extend it
to a class of doubling measures.

\begin{prop}\label{lemmaCZD-w} Let $n\ge 1$ and $1\le p< \infty$.
Let  $w\in L^1_{\rm loc}(\re^n)$, $w>0$ a.e., be such that
$d\mu=w\,dx$ is a Borel doubling measure \textup{(}here we do not
need that $w$ is a Muckenhoupt weight\textup{)}. Assume that the
measure $\mu$ supports an $L^p$ Poincar\'e inequality, that is,
\begin{equation}\label{Poincare-w}
\Big(\aver{B}|f-m_{B}f|^p \, d\mu\Big)^{\frac1p}
\le
C\, r(B)\,
\Big(\aver{B}|\nabla f|^p \, d\mu\Big)^{\frac1p}
\end{equation}
for all locally Lipschitz functions $f $   and all balls $B$ with
radius $r(B)$. Here $m_{B}f$ is the average of $f$ with respect to
$\mu$ on $B$. Assume that  $f\in \cals$ is  such that $\|\nabla
f\|_{L^p(\mu)} <\infty.$\footnote[2]{We avoid here regularity issues
by taking a smooth $f$.} Let $\alpha>0$. Then, one can find a
collection of balls $\{B_i\}_i$, smooth functions $\{b_i\}_i$ and a
function $g\in L^1_{\rm loc}(\RR^n,\mu)$ such that
\begin{equation}\label{eqcsds1}
f= g+\sum_i b_i
\end{equation} and the following properties hold:
\begin{equation}\label{eqcsds2}
|\nabla g(x)| \le C\alpha,
\quad
\text{for $\mu$-a.e. } x \ \footnote[3]{The gradient of $g$ exists
$\mu$-almost everywhere, that is almost everywhere for the Lebesgue
measure. In fact, a similar argument shows that $g$ is almost
everywhere equal to a Lipschitz function $\tilde g$. Hence, $\nabla
g$ coincide almost everywhere with the distributional gradient of
$\tilde g$.}
\end{equation}
\begin{equation}\label{eqcsds3}
\supp b_i \subset B_{i}
\quad\text{and}\quad
\int_{B_i} |\nabla b_i|^p\, d\mu \le C\alpha^p \mu(B_i),
\end{equation}
 \begin{equation}\label{eqcsds4}
\sum_i \mu(B_i) \le C\alpha^{-p} \int_{\re^n} |\nabla f|^p\, d\mu ,
\end{equation}
\begin{equation}\label{eqcsds5}
\sum_i \bigchi_{B_i} \le N,
\end{equation} where $C$ and  $N$ depends
only on dimension, the doubling constant of $\mu$ and $p$. Assuming
furthermore  that $\mu$ supports an $L^p-L^q$ Poincar\'e inequality
with $p\le q<\infty$, that is,
\begin{equation}
\Big(\aver{B}|f-m_{B}f|^q \, d\mu\Big)^{\frac1q}
\le
C\, r(B)\, \Big(\aver{B}|\nabla f|^p \, d\mu\Big)^{\frac1p}
\end{equation}
for all $f $ locally Lipschitz and all ball $B$. Then
\begin{equation}\label{eq20bis}
\Big(\aver{B_i} |b_i|^q\,d\mu\Big)^\frac1q
\lesssim
\alpha\, r(B_{i}).
\end{equation}

 \end{prop}

Since $A_{p}$ weights support an $L^p-L^q$ Poincar\'e inequality for some $q>p$,  the
latter result applies to any  $w\in A_{\infty}$ and $p > r_{w}$.

\section{Proofs of the main results}\label{section:proofs2}

We prove Theorem \ref{theor:B-K:small}, Theorem
\ref{theor:small:comm}, and Proposition \ref{lemmaCZD-w}.

\subsection{Proof of Theorem \ref{theor:B-K:small}}
\label{subsec:proof:theor:B-K:small}

We follow closely the proof in \cite{Aus-mem} (we include it since it
will be needed for the next section). By Marcinkiewicz interpolation
theorem, it suffices to show that $T$ is of weak-type $(p_0,p_0)$.
Let $f\in \D$ (so $f\in L^{p_0}(\mu)$) and $\alpha>0$. By the
Calder\'on-Zygmund decomposition (see \cite{CW} or
\cite{Ste}) for $|f|^{p_0}$ at height $\alpha^{p_0}$ it follows that
there exist a collection of balls $\{B_i\}_i$ and functions $g$,
$\{h_i\}_i$ such that $f=g+\sum_i h_i$ and the following properties
hold:
\begin{equation}\label{CZ:g}
\|g\|_{L^\infty(\mu)}
\le
C\,\alpha,
\end{equation}
\begin{equation}\label{CZ:hi}
\supp h_i\subset B_i,
\qquad
\Big(\aver{B_i} |h_i|^{p_0}\,d\mu\Big)^{\frac1{p_0}}
\le
C\,\alpha,
\end{equation}
\begin{equation}\label{CZ:Bi}
\sum_i \mu(B_i)
\le
C\, \alpha^{-p_0}\,\int_{\re^n}|f|^{p_0}\,d\mu,
\end{equation}
\begin{equation}\label{CZ:overlap}
\sum_i \bigchi_{B_i}\le N,
\end{equation}
where $C$ and $N$ depends on $\mu$, $n$ and $p_0$. We write
$r_i=r(B_i)$ and control $T f$ by
$$
|T f|
\le
|T g|+\Big|T\Big(\sum_i \A_{r_i}\,h_i\Big)\Big| +
\sum_i|T(I-\A_{r_i})h_i|
=
F_1+F_2+F_3.
$$
We  estimate $\mu\{F_{i}>\alpha/3\}$. For $F_1$, since $T$ is of
weak-type $(q_0,q_0)$ and \eqref{CZ:g}
\begin{equation} \label{CZ:F1}
\mu\{ F_1>\alpha/3\}
\lesssim
\frac1{\alpha^{q_0}}\int_{\re^n}|g|^{q_0}\,d\mu
\lesssim
\frac1{\alpha^{p_0}}\,\int_{\re^n}|g|^{p_0}\,d\mu
\lesssim
\frac1{\alpha^{p_0}}\,\int_{\re^n}|f|^{p_0}\,d\mu,
\end{equation}
where we have used that \eqref{CZ:overlap}, \eqref{CZ:hi},
\eqref{CZ:Bi} yield
$$
\int_{\re^n}\Big|\sum_i h_i\Big|^{p_0}\,d\mu
\lesssim
\sum_i \int_{B_i} |h_i|^{p_0}\,d\mu
\lesssim
\alpha^{p_0}\,\sum_i\mu(B_i)
\lesssim
\int_{\re^n}|f|^{p_0}\,d\mu.
$$
For  $F_{2}$, we first use that $T$ is of weak-type
$(q_0,q_0)$,
\begin{equation}\label{F2-1}
\mu\{F_{2}>\alpha/3\}
\lesssim
\frac1{\alpha^{q_0}}\,\int_{\re^n} \Big|\sum_i
\A_{r_i}h_i\Big|^{q_0}\,d\mu.
\end{equation}
To compute the $L^{q_0}$-norm we dualize against
$0\le u\in L^{q_0'}(\mu)$ with $\|u\|_{L^{q_0'}(\mu)}=1$. We use
\eqref{small:A}, \eqref{CZ:hi}, \eqref{CZ:overlap}
\begin{align}
\int_{\re^n}\Big|\sum_i \A_{r_i}h_i\Big|\,u\,d\mu
&\lesssim
\sum_i\sum_{j=1}^\infty 2^{j\,D}\,\mu(B_i)\,
\Big(\aver{C_j(B_i)} |\A_{r_i} h_i|^{q_0}\,d\mu\Big)^{\frac1{q_0}}\,
\Big(\aver{2^{j+1}\,B_i}u^{q_0'}\,d\mu\Big)^{\frac1{q_0'}}
\nonumber
\\
&\lesssim
\sum_i\sum_{j=1}^\infty 2^{j\,D}\,\mu(B_i)\, \alpha_j\,
\Big(\aver{B_i} |h_i|^{p_0}\,d\mu\Big)^{\frac1{p_0}}\,
\essinf_{y\in B_i} M_\mu\big(u^{q_0'}\big)^{\frac1{q_0'}}(y)
\nonumber
\\
&\lesssim
\alpha\, \int_{\re^n} \sum_i \bigchi_{B_i}
M_\mu\big(u^{q_0'}\big)^{\frac1{q_0'}}\,d\mu
\lesssim
\alpha\, \int_{\cup_i B_i}
M_\mu\big(u^{q_0'}\big)^{\frac1{q_0'}}\,d\mu \nonumber
\\
&\lesssim
\alpha\, \mu(\cup_i
B_i)^{\frac1{q_0}}\,\big\|u^{q_0'}\big\|_{L^1(\mu)}^{\frac1{q_0'}}
=
\alpha\, \mu(\cup_i B_i)^{\frac1{q_0}},\label{F2-2}
\end{align}
where we have used Kolmogorov's lemma and the weak-type $(1,1)$ for
the Hardy-Littlewood maximal function $M_\mu$ (this idea is borrowed
from
\cite{HM}). Next, we take the supremum on $u$ and plug the obtained
estimate into \eqref{F2-1}:
\begin{equation}\label{CZ:F2}
\mu\{F_{2}>\alpha/3\}
\lesssim
\mu(\cup_i B_i)
\lesssim
\frac1{\alpha^{p_0}}\,\int_{\re^n}|f|^{p_0}\,d\mu,
\end{equation}
where we have used \eqref{CZ:Bi}. Next, we consider $F_{3}$.  By
\eqref{small:T:I-A}, \eqref{CZ:hi} and \eqref{CZ:Bi}
\begin{align}
\lefteqn{\hskip-.7cm
\mu\big(( \re^n\setminus \cup_i 4\,B_i) \cap \{F_{3}>\alpha/3\}\big)
\le
\frac3{\alpha}\,\sum_i\int_{\re^n\setminus 4\,B_i}
|T(I-\A_{r_i})h_{i}|\,d\mu} \nonumber
\\
&\lesssim
\frac1{\alpha}\,\sum_i\sum_{j=2}^\infty 2^{j\,D}\,\mu(B_i)
\,\Big(\aver{C_j(B_i)}|T(I-\A_{r_i})h_{i}|\,d\mu\Big) \nonumber
\\
&\lesssim
\frac1{\alpha}\,\sum_i\sum_{j=2}^\infty
2^{j\,D}\,\mu(B_i)\,\alpha_j\,
\Big(\aver{B_i}|h_{i}|^{p_0}\,d\mu\Big)^{\frac1{p_0}}
\lesssim
\frac1{\alpha^{p_0}}\,\int_{\re^n}|f|^{p_0}\,d\mu. \label{CZ:F3}
\end{align}
Gathering \eqref{CZ:F1}, \eqref{CZ:F2}, \eqref{CZ:F3}, and using
\eqref{CZ:Bi} we conclude that
$$
\mu\{x\in \re^n:|T f(x)|>\alpha\}
\lesssim
\frac1{\alpha^{p_0}}\,\int_{\re^n}|f|^{p_0}\,d\mu.
$$

\subsection{Proof of Theorem \ref{theor:small:comm}}
\label{subsec:proof:theor:small:comm} The basic ingredient is the
following consequence of John-Nirenberg's inequality: for any ball
$B$, $0<s<\infty$ and $j\ge 0$,
\begin{equation}\label{JN:mu}
\Big(\aver{2^j\,B} |b-b_B|^{s}\,\,d\mu\Big)^{\frac1s}
\lesssim
(1+j)\, \|b\|_{\BMO(\mu)}.
\end{equation}

\begin{lemma}
Assume  \eqref{small:A} and \eqref{small:T:I-A:comm} of Theorem \ref{theor:small:comm}.  Let $p_{0}<p<q<q_{0}$. Let  $b\in L^\infty(\mu)$ with $\|b\|_{\BMO(\mu)}=1.$ Then for all ball  $B$ with radius $r$, all functions $f$ supported in $B$ and $m\in \NN$, $m\ge 1$,
\begin{equation}
\label{lemma:smallcomm1}
\Big( \aver{B} |(b-b_{4\,B})^m\, f|^{p_{0}}\,d\mu\Big)^{\frac1{p_{0}}}
\lesssim
\Big(\aver{B} |f|^{p}\,d\mu\Big)^{\frac1{p}},
\end{equation}
 for $j\ge 1$,
\begin{equation}
\label{lemma:smallcomm2}
\Big( \aver{C_{j}(B)} |(b-b_{4\,B})^m\, \A_{r} f|^{q}\,d\mu\Big)^{\frac1{q}}
\lesssim
j^m\, \alpha_{j}\,\Big(\aver{B} |f|^{p_0}\,d\mu\Big)^{\frac1{p_0}}
\end{equation}
and for  $j\ge 2$,
\begin{equation}\label{lemma:smallcomm3}
 \aver{C_j(B)} |(b-b_{4\,B})^m \, T(I-\A_{r})f|\,d\mu
\lesssim
j^m\, \alpha_j\,\Big(\aver{B} |f|^{p_0}\,d\mu\Big)^{\frac1{p_0}},
\end{equation}
where the constants involved are independent of $b$ and $f$.
\end{lemma}

The proof of \eqref{lemma:smallcomm1} is a direct application of
H\"older inequality and \eqref{JN:mu}. Next,  using that $q<q_{0}$,
\eqref{lemma:smallcomm2} follows from H\"older inequality,
\eqref{small:A}  and \eqref{JN:mu}.  Eventually,
\eqref{lemma:smallcomm3} is a consequence of H\"older inequality,
\eqref{small:T:I-A:comm}  as $r>1$  and \eqref{JN:mu}.

We begin the proof of Theorem \ref{theor:small:comm}.
As before it is enough to consider the case $b\in L^\infty(\mu)$
obtaining the desired estimates with a constant independent of $b$. Let us observe that here we assume
that $T$ is of weak-type $(q_0,q_0)$ in place of being bounded on
$L^{q_0}$. This changes slightly Lemma \ref{lemma:comm-apriori}.
Namely, in $(a)$ one  obtains that $T_b^kf\in L^{q_0,\infty}(\mu)$.
The proof of $(b)$ changes in the following way: one shows that
$T\big((b_N)^m\,f-b^m\,f\big)\longrightarrow 0$ in
$L^{q_0,\infty}(\mu)$ which also implies the convergence almost
everywhere for a subsequence. From here the proof can be carried out
in the same manner.

When $b\in L^\infty(\mu)$, all the formal computations below make
sense. Notice that by homogeneity, it suffices to consider the case
$\|b\|_{\BMO(\mu)}=1$. By Marcinkiewicz interpolation theorem, it
suffices to show that $T_b^k$ is of weak-type $(p,p)$ for all
$p_{0}<p<q_{0}$ because $T_{b}^k$ is sublinear. We proceed by
induction and assume that we have proved that $T_b^m$ is of weak-type
$(p,p)$ for all $p_{0}<p<q_{0}$ and $m=0, \ldots, k-1$, the case
$m=0$ being covered by Theorem \ref{theor:B-K:small}.

Fix $p$ so that $p_0<p<q_0$ and let $q$ with $p<q<q_{0}$.  Let $f\in
L^\infty_c$ (so $f\in L^{p}(\mu)$) and $\alpha>0$. By the
Calder\'on-Zygmund decomposition (see \cite{CW} or \cite{Ste}) for
$|f|^{p}$ at height $\alpha^{p}$ it follows that there exist a
collection of balls $\{B_i\}_i$, a collection  of functions
$\{h_i\}_i$ and a function $g$ such that $f=g+\sum_i h_i$ and
\eqref{CZ:g}, \eqref{CZ:hi}, \eqref{CZ:Bi} \eqref{CZ:overlap} hold
with $p$ in place of $p_0$. We wish to estimate $\mu\{ |T_{b}^k  f |>
\alpha\}.$ First, we have $$|T^k_{b}f| \le |T^k_{b}g| +
\Big|T^k_{b}\Big(\sum_{i} h_{i}\Big)\Big|.$$ By the weak-type
$(q_{0},q_{0})$ of $T_{b}^k$,
 \begin{equation} \label{CZ:F1:comm}
\mu\{|T^k_{b}g|>\alpha/2\}
\lesssim
\frac1{\alpha^{q_0}}\int_{\re^n}|g|^{q_0}\,d\mu
\lesssim
\frac1{\alpha^{p}}\,\int_{\re^n}|f|^{p}\,d\mu,
\end{equation}
where the last inequality follows as in \eqref{CZ:F1}. Next, set
$h_{i,b}^m=(b_{4\,B_i}-b)^m\, h_{i}$ and $r_{i}=r(B_{i})$. Then
\begin{align*}
\Big|T_b^k\Big(\sum_i \,h_i\Big)(x)\Big|
&
\le
\sum_{m=0}^k  C_{k,m}  \Big|T\Big(\sum_i (b(x) - b_{4\,
B_{i}})^{k-m}\,  \A_{r_{i}}h_{i,b}^m\Big)(x)\Big|
\\
&
\qquad +
\sum_{m=0}^k  C_{k,m} \sum_i |b(x) - b_{4\, B_{i}}|^{k-m}
\Big|T\Big( (I-\A_{r_{i}})h_{i,b}^m\Big)(x)\Big|
\end{align*}
The $m$-th term in the first sum is bounded by $\sum_{\ell=0}^{k-m}
c_{\ell}^m F_{m,\ell}(x)$ with
$$
F_{m,\ell}(x)= \Big|T_{b}^{k-m-\ell}\Big(\sum_i  (b-b_{4\, B_{i}})^\ell\, \A_{r_{i}}h_{i,b}^m\Big)(x)\Big|.
$$
Fix $\ell=m=0$ and $A$ some large number depending just on $k$.  Then
the estimate of $\mu\{F_{0,0}>\alpha/A\}$ is done as for the term
$F_{2}$ in the proof of Theorem \ref{theor:B-K:small}, using the
weak-type $(q_{0},q_{0})$ of $T_{b}^{k}$.  Next, fix $\ell, m$ with
$m+\ell>0$.  Then, the induction hypothesis implies that
$T_{b}^{k-m-\ell}$ is of weak-type $(q,q)$. Hence,  the estimate of
$\mu\{F_{m,\ell}>\alpha/A\}$ is done as for the term $F_{2}$ in the
proof of Theorem \ref{theor:B-K:small}, by replacing $q_{0}$ by $q$
and using \eqref{lemma:smallcomm2} with $f=h_{i,b}^m$ and then
\eqref{lemma:smallcomm1} with $f=h_{i}$.

It remains  to estimate  $\mu\{G_{m,\ell}>\alpha/A\}$ with
$$
G_{m,\ell}(x)= \sum_i  |b(x) - b_{4\, B_{i}}|^{k-m}   \Big|T\Big((I-\A_{r_{i}})h_{i,b}^m\Big)(x)\Big|.
$$
We proceed as for the term $F_{3}$ in the proof of Theorem
\ref{theor:B-K:small}, using \eqref{lemma:smallcomm3} with
$f=h_{i,b}^m$ and then \eqref{lemma:smallcomm1} with $f=h_{i}$. We
leave details to the reader.

\begin{remark}\label{remark:proof:multi-comm:p-small}\rm
The latter argument can be carried out for the multilinear
commutators introduced above. We give some of the ideas leaving
the precise computations to the reader. As before, it suffices to
consider the case $b_m\in L^\infty$ with $\|b_m\|_{\BMO(\mu)}=1$
for all $1\le m\le k$. Given $\sigma\subset\{1,\dots,k\}$, we
write $\pi_{i,\vec{b}_\sigma}=\prod_{j\in \sigma}
\big(b_j-(b_j)_{4\,B_i}\big)$ and $h_{i,\vec{b}_\sigma}=h_i\,\pi_{i,\vec{b}_\sigma}$
Here, when $\sigma=\emptyset$ we understand that
$\pi_{i,\vec{b}_\sigma}=1$ and $h_{i,\vec{b}_\sigma}=h_i$. Thus,
combining the preceding ideas with
\cite[p. 684]{PT} we have
$$
|T_{\vec{b}} f|
\le
|T_{\vec{b}} g|+ \sum_{\sigma_1,\sigma_2,\sigma_3}
\Big|T_{\vec{b}_{\sigma_1}}\Big(\sum_i \pi_{i,\vec{b}_{\sigma_2}}
\,\A_{r_i} h_{i,\vec{b}_{\sigma_3}}\Big)\Big| +
\sum_{\sigma_1,\sigma_2}\sum_i|\pi_{i,\vec{b}_{\sigma_1}}| \,
|T(I-\A_{r_i})h_{i,\sigma_2}|,
$$
where the first sum  (resp. the second sum) runs over all partitions
of $\{1,\dots,k\}$ in three (resp. two) pairwise disjoint sets
$\sigma_1,\sigma_2,\sigma_3$ (resp. $\sigma_1, \sigma_2$).

The estimate for the first term is obtained as in \eqref{CZ:F1:comm}.
The second term is treated as $F_{m,l}$ above (notice that the case
$\sigma_1=\{1,\dots,k\}$, $\sigma_2=\sigma_3=\emptyset$ is handled
differently as happened before). Finally, the third term is estimated
as $G_{m,l}$ above. Full details are left to the reader.
\end{remark}

\subsection{Proof of Proposition \ref{lemmaCZD-w}}

Let $\Omega= \{x \in \re^n: M_{\mu}(|\nabla f|^p)(x) >\alpha^p\}$
where $M_{\mu}$ is the uncentered maximal operator over
cubes\footnote[2]{We freely change balls to cubes.} of $\re^n$ with
respect to $\mu$. If $\Omega$ is empty, then set $g=f$. Otherwise,
since $\mu$ is doubling it follows that $M_\mu$ is of weak-type
$(p,p)$ and so
\begin{equation*}
|\Omega| \le C\alpha^{-p} \int_{\re^n} |\nabla f|^p\, d\mu.
\end{equation*} Let $F$ be the complement of $\Omega$. By the
Lebesgue differentiation theorem, $|\nabla f| \le \alpha$
$\mu$-almost everywhere on $F$.

\begin{lemma}\label{lemmalipschitzonF}
One can redefine $f$ on a $\mu$-null set of $F$ so that  for all $x
\in F$, and for all cube $Q$ centered at $x$,
\begin{equation}\label{eq19}
|f(x) - m_Qf| \le C\alpha \ell(Q) \end{equation} where $\ell(Q)$ is
the sidelength of $Q$. Furthermore, for all $x,y \in F$,
\begin{equation}\label{eq20}
|f(x) - f(y)| \le C\alpha |x-y|.
\end{equation}
The constant $C$ depends only on dimension, the doubling constant of
$\mu$ and $p$.
\end{lemma}

\begin{proof}[Proof of  Lemma \ref{lemmalipschitzonF}]
Let $x$ be a point in $F$. Fix a  cube $Q$ with center $x$ and let
$Q_k$ be co-centered cubes with $\ell(Q_k) =2^{-k}\,\ell(Q)$ for
$k\ge 1$. Then, by Poincar\'e's inequality
\begin{equation}\label{f-Qk}
|m_{Q_{k+1}}f - m_{Q_{k}}f|
\lesssim
\aver{Q_{k}} |f-m_{Q_{k}}f|\,d\mu
\lesssim
\ell(Q_{k}) \Big(\aver{Q_{k}} |\nabla f|^p d\mu\Big)^{\frac1p}
\lesssim
2^{-k}\ell(Q)\alpha
\end{equation}
since $x\in Q_{k}\cap F$. This easily implies that $\{m_{Q_k}
f\}_{k\ge 1}$ is a Cauchy sequence and so it converges as $k\to
\infty$ or what is the same as $\ell(Q_k)\to 0$. The Lebesgue
differentiation theorem implies that $m_{Q_k} f\longrightarrow f(x)$
whenever $x$ is a Lebesgue point of $f$, that is $\mu$-almost
everywhere.  If $x$ is not a Lebesgue point, it is easy to show that
$\lim m_{Q_k}f$ does not depend on $Q$ (the original cube).  Hence,
we redefine $f(x)$ as the value of this limit. With this new
definition, summing over $k\ge 1$ on \eqref{f-Qk} one gets
\eqref{eq19}.

To see \eqref{eq20}, let $x,y\in F$ and $Q_x$ be the cube centered at
$x$ with sidelength $2\,|x-y|$ and $Q_y$ be the cube centered at $y$
with sidelength $4\,|x-y|$. It is easy to see that $Q_x \subset Q_y$.
As in \eqref{f-Qk}, one can see that $|m_{Q_{x}}f - m_{Q_{y}}f|
\le C\alpha |x-y|$. Hence by the triangle inequality and
\eqref{eq19}, one obtains \eqref{eq20} readily.
\end{proof}

Let us continue the proof of Lemma \ref{lemmaCZD-w}. Let $\{Q_i\}_i$
be a Whitney decomposition of $\Omega$ by dyadic cubes. Hence,
$\Omega$ is the disjoint union of the $Q_i$'s, the cubes
$2\,Q_i\subset \Omega$ have bounded overlap, and the cubes $4\,Q_i$
intersect $F$. As usual, $\lambda\, Q$ is the cube co-centered with
$Q$ with sidelength $\ell(\lambda\,Q)=\lambda\,\ell(Q)$. Hence
\eqref{eqcsds4} and \eqref{eqcsds5} are satisfied by the cubes
$2\,Q_i$.

Let us now define the functions $b_i$ and show \eqref{eqcsds3}. Let
$\{\calX_i\}_i$ be a partition of unity on $\Omega$ associated to the
covering $\{Q_i\}_i$ so that for each $i$, $\calX_i$ is a $C^\infty$
function supported in $2\,Q_i$ with $\|\calX_i\|_\infty + \ell_i\,
\|\nabla \calX_i\|_\infty
\le c(n)$, $\ell_i$ being the sidelength of $Q_i$. Set
$$
b_i = (f-m_{2Q_{i}}f)\,\calX_i.
$$
It is clear that $b_i$ is supported in $2Q_i$. Since $\nabla
\big(\,(f-m_{2Q_{i}}f)\calX_i\big) = \calX_i \nabla f +
(f-m_{2Q_{i}}f)\nabla\calX_i$, we have  by the $L^p$ Poincar\'e
inequality, the fact that the average of $|\nabla f|^p$ on $4\,Q_i$
is controlled by $\alpha^p$ (since $4\,Q_i$ meets $F$) and the
doubling property that
$$\int_{2Q_i} |\nabla \big(\,(f-m_{2Q_{i}}f)\calX_i\big)|^p\, d\mu  \le
C\alpha^p \mu(2Q_i).$$ Thus \eqref{eqcsds3} is proved.

It remains to obtain \eqref{eqcsds1} and \eqref{eqcsds2}. To do so,
we introduce an auxiliary function  $h= \sum_i m_{2Q_{i}}f \
\nabla\calX_i$, for which we claim that $h \le C\alpha$ on  $\RR^n$.
First,  note that this sum is locally finite in $ \Omega$ and
vanishes on $F$, hence $h$ well-defined on $\RR^n$. Note also that
$\sum_i \calX_i$ is 1 on $\Omega$ and 0 on $F$. Since it is also
locally finite we have $\sum_i  \nabla\calX_i=0$ in $\Omega$.  Fix $x
\in \Omega$. Let $Q_j$ be the Whitney cube containing $x$ and let
$I_x$ be the set of indices $i$ such that $x \in 2\,Q_i$. We know
that $\# I_x \le N$. Also for $i \in I_x$ we have that
$C^{-1}\,\ell_i
\le  \ell_j \le C\,\ell_i$ where the constant $C$ depends only on
dimension (see \cite{Ste}).  We also have  $|m_{2Q_{i}}f-m_{2Q_{j}}f
|\le C\,\ell_j\,\alpha$ (embed $2\,Q_{i}$ and $2\,Q_{j}$ in some
dilate of $Q_{j}$ and apply Poincar\'e's inequality as in \eqref{f-Qk}
and the definition of $F$). Hence,
$$
|h(x)| = \left|\sum_{i \in I_x} (m_{2Q_{i}}f-m_{2Q_{j}}f)
\nabla\calX_i(x)\right| \le C \sum_{i \in I_x}
|m_{2Q_{i}}f-m_{2Q_{j}}f |\ell_i^{-1} \le C N \alpha.$$

We are ready to prove \eqref{eqcsds1} and \eqref{eqcsds2}. Set
$g=f-\sum b_{i}$.  This function is defined $\mu$-almost everywhere,
hence  \eqref{eqcsds1} trivially holds.  Next, we claim that $\nabla
g= {\bf 1}_F (\nabla f)  + h $ $\mu$-almost everywhere where ${\bf 1}_{E}$ is the indicator function of a set $E$. Admitting
this, for $\mu$-a.e.~$x\in F$, we have that $|\nabla g(x)|=|\nabla
f(x)|\le M_\mu
\big(|\nabla f|^p)(x)^p\le \alpha$, and for $\mu$-a.e.~$x\in
\Omega$,  $|\nabla g(x)|=| h(x)|\le C\, N\, \alpha$.  To  conclude
the proof of \eqref{eqcsds2}, it remains to see the claim. First,
observe that $\sum_i b_{i}$ converges in $L^p_{\rm loc}(\RR^n, \mu)$.
Indeed,  fix a compact set $K$ and observe that  the sidelengths of
the cubes $Q_{i}$ meeting $K$ are bounded. Since
$\|b_{i}\|_{L^p(\mu)}^p \le C\ell_{i}^p \alpha ^p \mu(Q_{i})$  and
$\sum_i \mu(Q_{i}) <\infty$, we obtain convergence of the series in
$L^p(K,\mu)$ from the bounded overlap property of the $Q_{i}$'s.
Next, it follows from \eqref{eqcsds3}, \eqref{eqcsds4} and
\eqref{eqcsds5}   that   $\sum_i |\nabla b_{i}|$ converges in
$L^p(\RR^n, \mu)$. We  invoke  \cite[Corollary 11]{FHK} (this is
where we use that $\mu$ is given by a weight) which implies that
$\nabla g$ exists almost everywhere (which is the same as
$\mu$-almost everywhere by the assumption on the weight) and is given
by $\nabla f - \sum_i \nabla b_{i}$. But as $\sum_i
\nabla\calX_i(x)=0$ for $x \in \Omega$, we have
$$
\nabla f = {\bf 1}_F (\nabla f) + {\bf 1}_\Omega(\nabla f)  = {\bf 1}_F  (\nabla
f)  + h + \sum_i \nabla b_i \quad \text{$\mu$-a.e.},
$$ and the claim follows.

It remains to prove \eqref{eq20bis} assuming  an  $L^p-L^q$  Poincar\'e
inequality. By the definition of  $b_{i}$ and similar computation as above,
$$
\Big(\aver{Q_i} |b_i|^q\,d\mu\Big)^\frac1q
\lesssim
\ell_i\, \Big(\aver{2\,Q_i} |\nabla f|^p\,d\mu\Big)^\frac1p
\lesssim
\ell_i\, \alpha.
$$

\end{document}